\documentclass[11pt]{amsart}
\usepackage{amsfonts}
\usepackage{amsmath}
\usepackage{amssymb}
\usepackage{tikz}
\usepackage{url,epsfig,eqngroup}
\usepackage{color}

\def\nat{\mathbb{N}}
\def\real{\mathbb{R}}
\def\cmplx{\mathbb{C}}
\def\ganz{\mathbb{Z}}

\def\eps{\varepsilon}

\def\dx{\mathrm{d}x}

\newcommand\cA{\mathcal{A}}

\newcommand\cF{\mathcal{F}}
\newcommand\cK{\mathcal{K}}
\newcommand\cL{\mathcal{L}}
\newcommand\cM{\mathcal{M}}
\newcommand\cN{\mathcal{N}}
\newcommand\cO{\mathcal{O}}

\newcommand\cR{\mathcal{R}}

\renewcommand{\Re}{\operatorname{Re}}
\renewcommand{\Im}{\operatorname{Im}}
 
\newcommand{\smalf}{\par\smallskip\noindent}
\newcommand{\medlf}{\par\medskip\noindent}
\newcommand{\biglf}{\par\bigskip\noindent}
 
\newtheorem{theorem}{Theorem}[section]
\newtheorem{definition}[theorem]{Definition} 
\newtheorem{lemma}[theorem]{Lemma}
\newtheorem{corollary}[theorem]{Corollary}

\newtheorem{remark}[theorem]{Remark}

\newtheorem{assumption}[theorem]{Assumption}

\begin{document}

\title{The PML-Method for a Scattering Problem for a Local Perturbation of an 
Open Periodic Waveguide}

\author{Andreas Kirsch and Ruming Zhang}

\footnotetext[1]{Department of
  Mathematics, Karlsruhe Institute of Technology (KIT), 76131 Karlsruhe, Germany. Email:
  andreas.kirsch@kit.edu. The author supported by the Deutsche Forschungsgemeinschaft (DFG, German Research Foundation) -- Project-ID 258734477 -- SFB 1173.}
\footnotetext[2]{Institute of Mathematics, Technische Universit{\"a}t Berlin, Stra\ss{}e des 17. Juni 136, 10623 Berlin, Germany. Email: zhang@math.tu-berlin.de. The author is supported by the Deutsche Forschungsgemeinschaft (DFG, German Research Foundation) Project-ID 433126998.}


\begin{abstract}
The perfectly matched layers method is a well known truncation technique for its efficiency and convenience in numerical implementations of wave scattering problems in unbounded domains.  In this paper, we study the convergence of the perfectly matched layers (PML) for wave scattering from a local perturbation of an open waveguide in $\real^2_+$, where the refractive index is a function  which is periodic along the axis of the waveguide and equals to one above a finite height. The problem is challenging due to the existence of guided waves, and a typical way to deal with the difficulty is to apply the limiting absorption principle. Based on the Floquet-Bloch transform and a curve deformation theory, the solution from the limiting absorption principle is rewritten as the integral of  a coupled family of quasi-periodic problems with respect to the quasi-periodicity parameter on a particularly designed curve. By comparing the Dirichlet-to-Neumann maps on a straight line above the locally perturbed periodic layer, we finally show that the PML method converges exponentially with respect to the PML parameter. Finally, the numerical examples are shown to illustrate the theoretical results.
\end{abstract}

\maketitle

\section{Introduction}
\label{s-intro}

Let $k>0$ be the wavenumber which is fixed throughout the paper and $n\in L^\infty(\real_+^2)$ 
be $2\pi-$periodic with respect to $x_1$ and equals to one for $x_2>h_0$ for some $h_0>0$. 
Furthermore, let $q\in L^\infty(\real^2_+)$ and $f\in L^2(\real_+^2)$ have compact supports in
$Q:=(0,2\pi)\times(0,h_0)$. It is the aim to solve
\begin{equation} \label{source-1}
\Delta u\ +\ k^2\,(n+q)\,u\ =\ -f\quad\text{in }\real_+^2\,,\quad u=0\mbox{ for }x_2=0\,,
\end{equation}
complemented by a suitable radiating condition stated below.
\smalf
The solution of (\ref{source-1}) is understood in the variational sense; that is,
\begin{equation} \label{var-Helmholtz}
\int\limits_{\real_+^2}\left[\nabla u\cdot\nabla\overline{\psi}-k^2\,(n+q)\,u\,\overline{\psi}
\right]\,dx\ =\ \int\limits_Qf\,\overline{\psi}\,dx
\end{equation}
for all $\psi\in H^1_0(\real_+^2)$ with compact support. By standard regularity theorems it is 
known that for $x_2>h_0$ the solution $u$ is a classical solution of the Helmholtz equation 
and thus analytic. 
\medlf
As mentioned above, a further condition is needed to assure uniqueness. As in 
\cite{Kirsc2017} we will derive the radiation condition by the limiting absorption 
principle; that is, the solution $u$ should be the limit (as $\eps>0$ tends to zero) of 
the solutions $u_\eps\in H^1(\real_+^2)$ corresponding to wave numbers $k+i\eps$ instead 
of $k$. 
\newline
In \cite{Hoang2011}, the author investigated the limiting absorption principle in deriving the physical solutions in periodic closed semi-waveguides. A further study of the structure of the solutions can be found in \cite{Fliss2015}.  The spectral decomposition of the propagating waves in closed periodic waveguides are given in \cite{Hohag2013,Zhang2019a}. Based on the singular perturbation theory, the radiation condition for wave propagation in closed periodic waveguides is proposed in \cite{Kirsc2017a}. Numerical methods are also developed based on the limiting absorption principle. One option is to approximate the Dirichlet-to-Neumann map on the periodicity boundaries by solving the cell problems. A Ricatti-equation based method was proposed for the first time in \cite{Joly2006} and then applied to other cases in \cite{Fliss2009a,Fliss2013,Fliss2021}; on the other hand, a doubling recursive process was also proposed in \cite{Yuan2007} and then followed by \cite{Ehrha2008a,Ehrha2009a,Ehrha2009b}.  For the closed waveguide scattering problems, there are only one type of singularity, i.e., the guided waves. The open waveguide problems are more complicated since the domain is unbounded in both directions and there are Rayleigh anomalies. Thus proper truncation techniques have to be proposed to truncate the infinite domain and deal with the singularities.
\newline
The situation is easier when the guided waves do not exist, for example, waves scattering from periodic surfaces.  With the help of the Floquet-Bloch transform, the problem is reduced into a coupled family of quasi-periodic problems, see \cite{Coatl2012,Hadda2016} for absorbing backgrounds and  see \cite{Lechl2015e,Lechl2016} for non-absorbing cases. Numerical solvers are also developed based on the Floquet-Bloch transform, see \cite{Lechl2016a,Lechl2017,Zhang2017e,Zhang2021b}. For more detailed explanation on the radiation conditions see \cite{Hu2021} and for the related boundary integral equation method see \cite{Yu2022}. For general cases when the guided waves exist, there are also works based on the singular perturbation theory, see \cite{Kirsc2017,Kirsc2019a,Kirsc2019b,Kirsc2022,Kirsc2023}. 
\newline
To simulate wave scattering from locally perturbed periodic layers, an efficient truncation technique is necessary to reduce the problem into a closed waveguide. In \cite{Lechl2016a,Lechl2017,Zhang2017e}, the exact but non-local transparent boundary condition was adopted and in \cite{Yu2022,Zhang2021b} the PML method was applied. The PML method, which was proposed by J.-P. Berenger in 1994 in \cite{Beren1994}, was well known for its convenience in numerical implementations. However, as it is not exact, the convergence is the key topic for this method. For  periodic surface scattering problems, the authors proved the exponential convergence for quasi-periodic problems except for some special quasi-periodicity in \cite{Chen2003}. For general rough surface scattering problems, the authors in \cite{Chand2009} proved the algebraic convergence but also asked a question that whether exponential convergence still holds in compact subsets. In \cite{Zhang2021b}, the author proved that the exponential convergence holds for all positive wavenumnbers except for half-integers. The higher-order algebraic convergence was proved for the exceptional cases in \cite{Zhang2023}. We will apply the PML method to truncate the open waveguide problem and study its convergence rate.
\newline
In this paper, we will adopt the curve deformation technique combined with the Floquet-Bloch transform, which was developed in \cite{Zhang2019b,Zhang2021a} to solve the closed waveguide scattering problems with existing guided waves, to rewrite the original solution by a simplified integral formulation with a couple family of quasi-periodic problems. On the specially designed curve, each quasi-periodic problem is uniquely solvable. The radiation condition, which was proposed in \cite{Kirsc2022}, can also be derived from the simplified integral formulation.  Then we apply the PML method to truncate the problem into a closed waveguide and the simplified integral formulation is again obtained for the PML problem. By comparing the Dirichlet-to-Neumann maps along the deformed curve, we finally prove the exponential convergence of the PML method.
\newline
The rest of the paper is organised as follows. In Section 2, we recalled the definitions and properties for the propagating modes (guided waves). In Section 3, the curve deformation technique combined with the Floquet-Bloch transform is applied to reformulate the scattering problems, and this approach is extended to locally perturbed cases in Section 4. In Section 5, the convergence analysis of the PML method is studied and the related numerical example are shown in Section 6.

\section{Propagating Modes}
\label{s-Prop}

First we recall some definitions. A function $\phi$ is $\alpha-$quasi-periodic (for some $\alpha\in\real$) if 
$\phi(x_1+2\pi)=e^{2\pi i\alpha}\phi(x_1)$ for all $x_1$. We define the layer 
$W_h:=\real\times(0,h)$, the space $H^1_{\alpha,0}(W_h):=\{u\in H^1_{loc}(W_h):u=0
\mbox{ for }x_2=0,\ u(\cdot,x_2)\mbox{ is $\alpha-$quasi-periodic for all }x_2\}$ of 
quasi-periodic functions, and the corresponding local space 
$H^1_{\alpha,0,loc}(\real_+^2)=\{u:u|_{W_h}\in H^1_{\alpha,0}(W_h)\mbox{ for all }h>0\}$.
\begin{definition} \label{d-exceptional}
(a) $\alpha\in\real$ is called a cut-off value if there exists $\ell\in\ganz$ with
$|\ell+\alpha|=k$.
\smalf
(b) $\alpha\in\real$ is called a propagative wave number (or quasi-momentum or Floquet 
spectral value) if there exists a non-trivial $\phi\in H^1_{\alpha,0,loc}(\real_+^2)$ such 
that 
\begin{equation} \label{mode}
\Delta\phi + k^2\,n\,\phi\ =\ 0\quad\text{in }\real_+^2
\end{equation}
satisfying the Rayleigh expansion
\begin{equation} \label{Rayleigh}
\phi(x)\ =\ \sum_{\ell\in\ganz}\phi_\ell(h_0)\,e^{i(\ell+\alpha)x_1+
i\sqrt{k^2-(\ell+\alpha)^2}x_2}\quad\mbox{for }x_2>h_0\,.
\end{equation}
Here, $\phi_\ell(h_0)=\frac{1}{\sqrt{2\pi}}\int_0^{2\pi}\phi(x_1,h_0)
e^{-i(\ell+\alpha)x_1}dx_1$ are the Fourier coefficients of $\phi(\cdot,h_0)$. The 
convergence is uniform for $x_2\geq h_0+\delta$ for all $\delta>0$. The functions $\phi$ 
are called propagating (or guided) modes. The branch of the square root is taken such 
that the square root is holomorphic in $\cmplx\setminus i\real_{\leq 0}$.
\end{definition}
If we decompose $k$ into $k=\hat{\ell}+\kappa$ with $\hat{\ell}\in\nat\cup\{0\}$ and 
$\kappa\in(-1/2,1/2]$ we observe that the cut-off values are given by $\pm\kappa+\ell$ for 
any $\ell\in\ganz$.
\smalf
Since with $\alpha$ also $\alpha+\ell$ for every $\ell\in\ganz$ is a propagative wave 
number we can restrict ourselves to propagative wave numbers in $(-1/2,1/2]$.
\newline 
We define the spaces of periodic functions with boundary conditions for $x_2=0$ by 
$H^1_{per,0}(Q):=\{\psi\in H^1(Q):\psi=0\mbox{ for }x_2=0,\ \psi\mbox{ is $2\pi-$periodic 
with respect to }x_1\}$ and, analogously, $H^1_{per,0}(Q^\infty)$ for 
$Q^\infty:=(0,2\pi)\times(0,\infty)$.
\newline 
For the proof that there exists at most a finite number of propagative wave number 
we need the following result.
\begin{lemma} \label{l-equiv} \mbox{ }
\begin{itemize}
\item[(a)] Let $u\in H^1_{\alpha,0,loc}(\real_+^2)$ be a $\alpha-$quasi-periodic solution of
$\Delta u+k^2nu=-f$ in $\real^2_+$ satisfying the Rayleigh expansion (\ref{Rayleigh}). Then 
$\tilde{u}(x):=e^{-i\alpha x_1}u(x)$ for $x\in Q$ is in $H^1_{per,0}(Q)$ and satisfies
\begin{eqnarray} 
& & \int\limits_Q\biggl[\nabla\tilde{u}\cdot\nabla\overline{\psi}-
2i\alpha\,\partial_1\tilde{u}\,\overline{\psi}-(k^2\,n-\alpha^2)\,\tilde{u}\,
\overline{\psi}\biggr]\,dx\ - \label{var-1} \\
& &  -\ i\sum_{\ell\in\ganz}\sqrt{k^2-(\ell+\alpha)^2}\,
\tilde{u}_\ell(h_0)\,\overline{\psi_\ell(h_0)} \nonumber \\ 
& = & \int\limits_Qe^{-i\alpha x_1}\,f(x)\,\overline{\psi(x)}\,dx
\quad\mbox{for all }\psi\in H^1_{per,0}(Q)\,. \nonumber
\end{eqnarray}
\item[(b)] If $\tilde{u}\in H^1_{per,0}(Q)$ solves (\ref{var-1}) then $u(x)=e^{i\alpha x_1}
\tilde{u}(x)$ for $x\in Q$ and its extension by the Rayleigh expansion (\ref{Rayleigh}) is 
in $H^1_{\alpha,0,loc}(\real_+^2)$ and satisfies (\ref{source-1}).
\item[(c)] The variational equation (\ref{var-1}) can be written as $(I-K_{k,\alpha})\tilde{u}=
r_\alpha$ in $H^1_{per,0}(Q)$ where $r_\alpha\in H^1_{per,0}(Q)$ and $K_{k,\alpha}$ is a 
compact operator from $H^1_{per,0}(Q)$ into itself. The operator depends continuously on 
$\alpha\in[-1/2,1/2]$ and $k>0$. Furthermore, for every $\hat{k}>0$ and $\hat{\alpha}\in
[-1/2,1/2]$ which is not a cut-off value with respect to $\hat{k}$ there exist neighborhoods 
$U,V\subset\cmplx$ of $\hat{k}$ and $\hat{\alpha}$, respectively, such that $K_{k,\alpha}$ 
depends analytically on $(k,\alpha)\in U\times V$. Finally, $r_\alpha$ depends analytically 
on $\alpha\in\cmplx$.
\end{itemize}
\end{lemma}
\textbf{Proof:} The parts (a) and (b) are standard. For (c) we choose
\begin{equation} \label{ast}
(u,v)_\ast\ :=\ \int\limits_Q\bigl[\nabla u\cdot\nabla\overline{v}+u\,\overline{v}
\bigr]\,dx\ +\ \sum_{\ell\in\ganz}|\ell|\,u_\ell(h_0)\,\overline{v_\ell(h_0)}
\end{equation}
as the inner product in $H^1_{per,0}(Q)$. The Theorem of Riesz implies the existence of 
$r\in H^1_{per,0}(Q)$ and an operator $K_{k,\alpha}$ from $H^1_{per,0}(Q)$ into itself that 
(\ref{var-1}) can be written as $(I-K_{k,\alpha})\tilde{u}=r_\alpha$. The compact imbedding of 
$H^1_{per,0}(Q)$ in $L^2(Q)$ and the fact that $\sqrt{k^2-(\ell+\alpha)^2}-|\ell|$ is 
bounded yields compactness of $K_{k,\alpha}$.
\newline
The fact that $\hat{\alpha}$ is not a cutoff value implies the existence of $c>0$ with 
$\left|\hat{k}^2-(\ell+\hat{\alpha})^2\right|\geq c$ for all $\ell\in\ganz$. We show that there 
exists $c_1>0$ and $\delta>0$ such that also $\left|\Re\left[k^2-(\ell+\alpha)^2\right]
\right|\geq c_1$ for all $\ell\in\ganz$ and all $k,\alpha\in\cmplx$ with $\left|k-\hat{k}\right|<\delta$ 
and $\left|\alpha-\hat{\alpha}\right|<\delta$. We consider three cases:
\begin{itemize}
\item[(i)] $|\ell|\geq\hat{k}+2$. Then $\Re[(\ell+\alpha)^2-k^2]=(\ell+\Re\alpha)^2-
(\Im\alpha)^2-(\Re k)^2+(\Im k)^2\geq(|\ell|-1)^2-(\Re k)^2-(\Im\alpha)^2\geq 
(\hat{k}+1)^2-(\Re k)^2-(\Im\alpha)^2\geq c_1$ for sufficiently small $\left|\Re k-\hat{k}\right|$ 
and $|\Im\alpha|$.
\item[(ii)] $|\ell|\leq\hat{k}+2$ and $\hat{k}^2-(\ell+\hat{\alpha})^2\geq c$. Then 
$\Re[k^2-(\ell+\alpha)^2]=(\Re k)^2-(\Im k)^2-(\ell+\Re\alpha)^2+(\Im\alpha)^2\geq
c+[(\Re k)^2-\hat{k}^2]-(\Im k)^2+[(\ell+\hat{\alpha})^2-(\ell+\Re\alpha)^2]\geq 
c+[(\Re k)^2-\hat{k}^2]-(\Im k)^2-|\hat{\alpha}-\Re\alpha|\cdot 2(\hat{k}+2)|
\hat{\alpha}+\Re\alpha|\geq c_1$ for sufficiently small $|k-\hat{k}|$ and 
$|\hat{\alpha}-\alpha|$.
\item[(iii)] $|\ell|\leq\hat{k}+2$ and $(\ell+\hat{\alpha})^2-\hat{k}^2\geq c$. Then
$\Re[(\ell+\alpha)^2-k^2]=(\ell+\Re\alpha)^2-(\Im\alpha)^2-(\Re k)^2+(\Im k)^2\geq
c+[\hat{k}^2-(\Re k)^2]-(\Im\alpha)^2 +[(\ell+\Re\alpha)^2-(\ell+\hat{\alpha})^2]\geq
c+[\hat{k}^2-(\Re k)^2]-(\Im\alpha)^2-|\hat{\alpha}-\Re\alpha|\cdot 2(\hat{k}+2)|
\hat{\alpha}+\Re\alpha|\geq c_1$ for sufficiently small $|k-\hat{k}|$ and 
$|\hat{\alpha}-\alpha|$.
\end{itemize}
Recalling our choice of the branch of the square root function we observe that the square 
roots in (\ref{var-1}) are holomorphic with respect to $k$ and $\alpha$ and their derivatives 
with respect to $k$ or $\alpha$ are bounded with respect to $\ell$. Therefore, the operator 
$K_{k,\alpha}$ depends analytically on $k$ and $\alpha$ in neighborhoods $U,V\subset\cmplx$ of
$\hat{k}$ and $\hat{\alpha}$, respectively. \qed
\biglf
Under the following assumption it can easily be shown as in, e.g., \cite{Kirsc2017} 
that every propagating mode $\phi$ corresponding to some propagative wave number 
$\alpha$ is evanescent; that is, $\phi_\ell(\pm h_0)=0$ for all $|\ell+\alpha|\leq k$; that 
is, there exist $c,\delta>0$ with $|\phi(x)|\leq c\,e^{-\delta x_2}$ for all $x_2>h_0$.
\begin{assumption} \label{assump1}
Let $|\ell+\alpha|\not=k$ for all propagative wave numbers $\alpha$ and all 
$\ell\in\ganz$; that is, the cut-off values are no propagative wave numbers.
\end{assumption}
\begin{lemma}
Under Assumption~\ref{assump1} there exists at most a finite number of propagative wave 
numbers in $[-1/2,1/2]$. Furthermore, if $\alpha$ is a propagative wave number with mode 
$\phi$ then $-\alpha$ is a propagative wave number with mode $\overline{\phi}$. Therefore, 
we can numerate the propagative wave numbers in $[-1/2,1/2]$ such they are given by 
$\{\hat{\alpha}_j:j\in J\}$ where $J\subset\ganz$ is symmetric with respect to $0$ and 
$\hat{\alpha}_{-j}=-\hat{\alpha}_j$ for $j\in J$. Furthermore, it is known that every 
eigenspace 
\begin{equation} \label{X_j}
\hat{X}_j\ :=\ \bigl\{ \phi\in H^1_{\hat{\alpha}_j,0,loc}(\real_+^2):\phi
\mbox{ satisfies (\ref{mode}) and (\ref{Rayleigh})} \bigr\}
\end{equation}
is finite dimensional with some dimension $m_j>0$. 
\end{lemma}
\textbf{Proof:} We recall that $\pm\kappa$ are the cut-off values in $[-1/2,1/2]$. We can 
cover the set $[-1/2,1/2]\setminus\{\kappa,-\kappa\}$ by at most three open sets $U_j\subset
\cmplx$ such that the operator $K_{k,\alpha}$ depends analytically on $\alpha\in\cup_jU_j$. 
Assume that there exists an infinite number of propagative wave numbers in $[-1/2,1/2]$. Then 
a subsequence converges to some $\alpha\in[-1/2,1/2]$ and, without loss of generality, this 
subsequence lies in one of the $U_j$. By the analytic Fredholm theory (see, e.g. 
\cite{Colto2019}) either every point of $U_j$ is a propagative wave number or the discrete 
set of propagative wave numbers has no limits point in $U_j$; that is, $\alpha\in 
[-1/2,1/2]\cap\partial U_j$ which implies that $\alpha$ coincides with $\kappa$ or 
$-\kappa$. In any case, by a perturbation argument, $\alpha$ is a propagative wave number 
which contradics Assumption~\ref{assump1}. \qed 
\biglf
In \cite{Kirsc2017} we have seen that a suitable basis of $\hat{X}_j$ is given by the 
eigenfunctions $\{\hat{\phi}_{\ell,j}:\ell=1,\ldots,m_j\}$ of the following selfadjoint 
eigenvalue problem
$$ -2i\int\limits_{Q^\infty}\partial_1\hat{\phi}\,\overline{\psi}\,dx\ =\ 
\lambda\,2k\int\limits_{Q^\infty}n\,\hat{\phi}\,\overline{\psi}\,dx
\quad\mbox{for all }\psi\in\hat{X}_j\,. $$
Therefore,
\begin{eqngroup}\begin{equation} \label{evp:a}
-2i\int\limits_{Q^\infty}\overline{\psi}\,\partial_1\hat{\phi}_{\ell,j}\,dx\ =\ 
\lambda_{\ell,j}\,2k\int\limits_{Q^\infty}n\,\hat{\phi}_{\ell,j}\,\overline{\psi}\,dx\,,\quad
\psi\in\hat{X}_j\,,
\end{equation}
for $\ell=1,\ldots,m_j$ with normalization
\begin{equation} \label{evp:b}
2k\int\limits_{Q^\infty}n\,\hat{\phi}_{\ell,j}\,\overline{\hat{\phi}_{\ell^\prime,j}}\,dx\ 
=\ \delta_{\ell,\ell^\prime}\,.
\end{equation}\end{eqngroup}

\section{The Limiting Absorption Solution for the Unperturbed Case}

First we study the unperturbed problem with a wavenumber of the form $k+i\eps$ for some 
$k>0$ and $\eps>0$. Then the Theorem of Lax-Milgram implies that the problem
$\Delta u_\eps+(k+i\eps)^2u_\eps=-f$ has a unique solution $u_\eps\in H^1(\real^2_+)$ with 
$u_\eps=0$ for $x_2=0$. We apply the (periodic) Floquet-Bloch transform to $u_\eps$. 
Therefore, 
$$ \tilde{u}(x,\eps,\alpha)\ :=\ \sum_{\ell\in\ganz}u_\eps(x_1+2\pi\ell,x_2)\,
e^{-i(x_1+2\pi\ell)\alpha} $$ 
is in $H^1_{per,0}(Q^\infty)$ and is the unique periodic (wrt $x_1$) solution of 
$\Delta\tilde{u}(\cdot,\eps,\alpha)+2i\alpha\,\partial_1\tilde{u}(\cdot,\eps,\alpha)+
[(k+i\eps)^2-\alpha^2]\tilde{u}(\cdot,\eps,\alpha)=-\tilde{f}(\cdot,\alpha)$ in 
$Q^\infty$ with $\tilde{u}(\cdot,\eps,\alpha)=0$ for $x_2=0$. 
We note that the function $\tilde{f}$ has the form $\tilde{f}(x,\alpha)=e^{-i\alpha x_1}f(x)$ 
because the support of $f$ is in $Q$. By the analog of Lemma~\ref{l-equiv} (equation 
(\ref{var-1})) the function $\tilde{u}$ satisfies
\begin{eqnarray} 
& & \int\limits_Q\biggl[\nabla\tilde{u}\cdot\nabla\overline{\psi}-
2i\alpha\,\partial_1\tilde{u}\,\overline{\psi}-((k+i\eps)^2\,n-\alpha^2)\,\tilde{u}\,
\overline{\psi}\biggr]\,dx\ - \label{var-2} \\
& &  -\ i\sum_{\ell\in\ganz}\sqrt{(k+i\eps)^2-(\ell+\alpha)^2}\,
\tilde{u}_\ell(h_0)\,\overline{\psi_\ell(h_0)} \nonumber \\ 
& = & \int\limits_Qf(x)\,e^{-i\alpha x_1}\,\overline{\psi(x)}\,dx\quad\mbox{for all }
\psi\in H^1_{per,0}(Q)\,. \nonumber
\end{eqnarray}
We write this in the form 
\begin{equation} \label{eq-operator}
(I-K_{\eps,\alpha})\tilde{u}(\eps,\alpha)\ =\ R_\alpha f\quad\mbox{in }H^1_{per,0}(Q)
\end{equation}
where
\begin{eqngroup}\begin{eqnarray}
\bigl((I-K_{\eps,\alpha})v,\psi\bigr)_\ast & = & \int\limits_Q\biggl[\nabla v\cdot
\nabla\overline{\psi}-2i\alpha\,\partial_1v\,\overline{\psi}-((k+i\eps)^2\,n-\alpha^2)\,v\,
\overline{\psi}\biggr]\,dx \label{KR:a} \\
& &  -\ i\sum_{\ell\in\ganz}\sqrt{(k+i\eps)^2-(\ell+\alpha)^2}\,
v_\ell(h_0)\,\overline{\psi_\ell(h_0)}\,, \nonumber \\
\bigl(R_\alpha g,\psi\bigr)_\ast & = & \int\limits_Qg(x)\,e^{-i\alpha x_1}\,
\overline{\psi(x)}\,dx \label{KR:b}
\end{eqnarray}\end{eqngroup}
for all $v,\psi\in H^1_{per,0}(Q)$ and $g\in L^2(Q)$.
\medlf
It is the aim to apply the following abstract representation theorem. First, we introduce the 
(punctured) cylinder $B_\delta:=\{(\eps,\alpha)\in[0,\delta)\times\cmplx:|\alpha|<\delta,\ 
\eps+|\alpha|>0\}$ in $\real\times\cmplx$.
\begin{theorem} \label{t-absorp-abstract}
Let $V\subset\cmplx$ be an open set containing $0$. Let $K(\eps,\alpha):H\to H$ be a family 
of compact operators from a (complex) Hilbert space $H$ into itself and 
$r(\eps,\alpha)\in H$ such that $(\eps,\alpha)\mapsto K(\eps,\alpha)$ and 
$(\eps,\alpha)\mapsto r(\eps,\alpha)$ are continuously differentiable on $[0,\eps_0]\times V$
for some $\eps_0>0$.\footnote{That is, the partial derivatives with respect to the real 
variable $\eps$ and the complex variable $\alpha$ exist in $(0,\eps_0)$ and $V$, 
respectively, and can be continued continuously into $[0,\eps_0]$ and $\overline{V}$, 
respectively.} Set $L(\eps,\alpha):=I-K(\eps,\alpha)$ and assume the following:
\begin{itemize}
\item[(i)] The null space $\cN:=\cN\bigl(L(0,0)\bigr)$ is not trivial and the 
Riesz number of of $L(0,0)$ is one; that is, the algebraic and geometric
multiplicities of the eigenvalue $1$ of $K(0,0)$ coincide. Let 
$P:H\to\cN\subset H$ be the projection operator onto $\cN$ corresponding to the 
direct decomposition $H=\cN\oplus\cR\bigl(L(0,0)\bigr)$, 
\item[(ii)] $A:=iP\frac{\partial}{\partial\eps}L(0,0)|_{\cN}:\cN\to\cN$ is 
selfadjoint and positive definite and $B:=-P\frac{\partial}{\partial\alpha}
L(0,0)|_{\cN}:\cN\to\cN$ is selfadjoint and one-to-one. 
\end{itemize}
Let $\bigl\{\lambda_\ell,\phi_\ell:\ell=1,\ldots,m\bigr\}$ be an 
orthonormal eigensystem of the following generalized selfadjoint problem in the 
$m-$dimensional space $\cN$:
\begin{equation} \label{ewp}
-B\phi_\ell\ =\ \lambda_\ell\,A\phi_\ell\quad\mbox{in $\cN$\quad with normalization}\quad 
\bigl(A\phi_\ell,\phi_{\ell^\prime}\bigr)_H\ =\ \delta_{\ell,\ell^\prime}\,.
\end{equation}
We set $M:=\{(\eps,\alpha)\in\real\times\cmplx:\sigma\Im\alpha\leq 0\}$ if all $\lambda_\ell$ 
have the same sign $\sigma\in\{+,-\}$, otherwise we set $M:=\real\times\real$. Then there 
exists $\delta>0$ such that:
\begin{itemize}
\item[(a)] For $(\eps,\alpha)\in M\cap B_\delta$ the equation $L(\eps,\alpha)u(\eps,\alpha)=
r(\eps,\alpha)$ has a unique solution $u(\eps,\alpha)\in H$, and $u(\eps,\alpha)$ has the 
form
\begin{equation} \label{rep-abstract}
u(\eps,\alpha)\ =\ u^b(\eps,\alpha)\ -\ \sum_{\ell=1}^m\frac{r_\ell}
{i\eps-\lambda_\ell\alpha}\,\phi_\ell\quad\mbox{for }(\eps,\alpha)\in M\cap B_\delta\,,
\end{equation}
where $r_\ell=\bigl(Pr(0,0),\phi_\ell\bigr)_H$ are the expansion coefficients of 
$A^{-1}Pr(0,0)$ with respect to the inner product $(A\cdot,\cdot)_H$; that is, 
$Pr(0,0)=\sum_{\ell=1}^mr_\ell\,A\phi_\ell$. Furthermore, $u^b(\eps,\alpha)$ depends 
continuously on $(\eps,\alpha)\in M\cap B_\delta$ and there exists $c>0$ with 
$\Vert u^b(\eps,\alpha)\Vert_H\leq c\Vert r\Vert_{C^1(\overline{B_\delta};H)}$ for all 
$(\eps,\alpha)\in M\cap B_\delta$.
\footnote{Here, $\Vert r\Vert_{C^1(\overline{B_\delta};H)}=\sup\limits_{(\eps,\delta)\in 
B_\delta}\Vert r(\eps,\delta)\Vert_H+\sup\limits_{(\eps,\delta)\in B_\delta}
\Vert\partial_\eps r(\eps,\delta)\Vert_H+\sup\limits_{(\eps,\delta)\in B_\delta}
\Vert\partial_\alpha r(\eps,\delta)\Vert_H$.}
\item[(b)] If all $\lambda_\ell$ have the same sign $\sigma\in\{+,-\}$ then $u^b(\eps,\cdot)$ 
depends analytically\footnote{Since the notions of a complex differentiable, a 
holomorphic, or an analytic function from a domain in $\cmplx$ into a Banach space coincide 
(see, e.g. \cite{}), we use the notion of analyticity in the following.} on 
$\alpha\in\bigl\{\alpha\in\cmplx:|\alpha|<\delta,\ \sigma\Im\alpha<0\bigr\}$ for every 
$\eps\in(0,\delta)$. 
\item[(c)] For $\eps=0$ the part $u^b(o,\cdot)$ depends analytically on $\alpha\in
\{\alpha\in\cmplx:|\alpha|<\delta\}$.
\end{itemize}
\end{theorem}
\textbf{Proof:} We project the equation $u(\eps,\alpha)-K(\eps,\alpha)u(\eps,\alpha)=
r(\eps,\alpha)$ onto $\cN$ and $\cR:=\cR(\bigl(L(0,0)\bigr)$. Decomposing $u(\eps,\alpha)$ 
into $u(\eps,\alpha)=u^N(\eps,\alpha)+u^R(\eps,\alpha)$ with $u^N(\eps,\alpha)\in\cN$ and 
$u^R(\eps,\alpha)\in\cR$ the equation is equivalent to the system
\begin{eqnarray*}
u^N(\eps,\alpha)-PK(\eps,\alpha)[u^N(\eps,\alpha)+u^R(\eps,\alpha)] & = &
Pr(\eps,\alpha)\,, \\
u^R(\eps,\alpha)-QK(\eps,\alpha)[u^N(\eps,\alpha)+u^R(\eps,\alpha)] & = &
Qr(\eps,\alpha)\,.
\end{eqnarray*}
The operator $[I-QK(0,0)]|_{\cR}:\cR\to\cR$ is invertible. Therefore, there exists 
$\delta_1>0$ such that $[I-QK(\eps,\alpha)]|_{\cR}:\cR\to\cR$ is invertible for all 
$(\eps,\alpha)\in B_\delta$. Set $T(\eps,\alpha):=[I-QK(\eps,\alpha)]|_{\cR}^{-1}:\cR\to\cR$.
Then we solve the second equation for $u^R(\eps,\alpha)$ and substitute it into the first 
equation which gives
$$ u^N(\eps,\alpha)-PK(\eps,\alpha)[I+T(\eps,\alpha)QK(\eps,\alpha]u^N(\eps,\alpha)\ =\
Pr(\eps,\alpha)+PK(\eps,\alpha)T(\eps,\alpha)Qr(\eps,\alpha) $$ 
in the finite dimensional space $\cN$. We abbreviate this as $C(\eps,\alpha)u^N(\eps,\alpha)
= g(\eps,\alpha)$. We note that $C(0,0)=0$ and $g(0,0)=Pr(0,0)$ and $\partial_\eps C(0,0)=
-P\partial_\eps K(0,0)|_{\cN}=-iA$ and $\partial_\alpha C(0,0)=
-P\partial_\alpha K(0,0)|_{\cN}=-B$.
\newline
We compare the equation $C(\eps,\alpha)u^N(\eps,\alpha)=g(\eps,\alpha)$ with the linearized 
equation
$$ [i\eps A+\alpha B]\tilde{u}(\eps,\alpha)\ =\ -g(0,0)\ =\ -Pr(0,0)\quad\mbox{in }\cN\,. $$
This equation is explicitly solved by 
$$ \tilde{u}(\eps,\alpha)\ =\ -\sum_{\ell=1}^m\frac{r_\ell}{i\eps-\lambda_\ell\alpha}\,
\phi_\ell\quad\mbox{for all }(\eps,\alpha)\in M\quad\mbox{where}\quad Pr(0,0)=
\sum_{\ell=1}^m r_\ell\,A\phi_\ell\,. $$
Indeed, if $M=\real\times\real$ then $|i\eps-\lambda_\ell\alpha|^2=\eps^2+
\lambda_\ell^2\alpha^2\geq\gamma(\eps^2+\alpha^2)$ for all $(\eps,\alpha)\in M$ with 
$\gamma=\min\{1,\lambda_1^2,\ldots,\lambda_m^2\}$. If all $\lambda_\ell$ have the same sign 
$\sigma\in\{+,-\}$ then $|i\eps-\lambda_\ell\alpha|^2=\eps^2+\lambda_\ell^2|\alpha|^2-
2\eps\lambda_\ell\Im\alpha\geq\gamma(\eps^2+|\alpha|^2)$ for all 
$(\eps,\alpha)\in\real\times\cmplx$ with $\sigma\Im\alpha\leq 0$.
\newline
Next we show the existence of $c>1$ such that
\begin{equation} \label{aux3}
\frac{1}{c}\,\Vert v\Vert_H\ \leq\ \sqrt{\eps^2+|\alpha|^2}\,
\bigl\Vert[i\eps A+\alpha B]^{-1}v\bigr\Vert_H\ \leq\ c\,\Vert v\Vert_H
\end{equation}
for all $v\in\cN$ and $(\eps,\alpha)\in M$. Indeed, let $u=[i\eps A+\alpha B]^{-1}v$. 
Then, as before, $u=\sum_{\ell=1}^m\frac{v_\ell}{i\eps-\lambda_\ell\alpha}\phi_\ell$ where 
$v_\ell$ are the coefficients in the expansion $A^{-1}v=\sum_{\ell=1}^mv_\ell\phi_\ell$. By 
the orthonormality of $\{\phi_\ell:\ell=1,\ldots,m\}$ with respect to $(Au,v)_H$ we have 
$(Au,u)_H=\sum_{\ell=1}^m\frac{|v_\ell|^2}{|i\eps-\lambda_\ell\alpha|^2}\leq
\frac{1}{\gamma(\eps^2+|\alpha|^2)}\sum_{\ell=1}^m|v_\ell|^2=\frac{1}{\gamma(\eps^2+
|\alpha|^2)}\,(v,A^{-1}v)_H$. Also, the norms $\sqrt{(Av,v)_H}$ and $\sqrt{(u,A^{-1}u)_H}$ 
are equivalent to $\Vert v\Vert_H$ and $\Vert u\Vert_H$, respectively, which proves 
(\ref{aux3}). 
\newline 
In particular we have that $\Vert\tilde{u}(\eps,\alpha)\Vert_H\leq
\frac{c}{\sqrt{\eps^2+|\alpha|^2}}$. We set $S_{\eps,\alpha}:=[i\eps A+\alpha B]^{-1}$ and 
note that $(\eps,\alpha)\mapsto S_{\eps,\alpha}$ is continuous on $M$ and 
$\Vert S_{\eps,\alpha}\Vert\leq\frac{c}{\sqrt{\eps^2+|\alpha|^2}}$ for all $(\eps,\alpha)\in 
M$. If all $\lambda_\ell$ have the same sign then $S_{\eps,\alpha}$ depends analytically on 
$\alpha$.
\smalf
Now we consider the difference $v=\tilde{u}(\eps,\alpha)-u^N(\eps,\alpha)$ and have
\begin{eqnarray*}
\bigl[i\eps A+\alpha B\bigr]v(\eps,\alpha) & = & g(\eps,\alpha)\ -\ g(0,0)\ +\
\bigl[C(\eps,\alpha)+i\eps A+\alpha B\bigr]v(\eps,\alpha) \\
& & -\ \bigl[C(\eps,\alpha)+i\eps A+\alpha B\bigr]\tilde{u}(\eps,\alpha)
\end{eqnarray*}
which we write in the form
\begin{eqnarray*}
& & \bigl(I-S_{\eps,\alpha}[C(\eps,\alpha)+i\eps A+\alpha B]\bigr)v(\eps,\alpha) \\
& = & S_{\eps,\alpha}[g(\eps,\alpha)-g(0,0)]\ -\
S_{\eps,\alpha}[C(\eps,\alpha)+i\eps A+\alpha B]\tilde{u}(\eps,\alpha) 
\end{eqnarray*}
From $\Vert S_{\eps,\alpha}\Vert\leq\frac{c}{\sqrt{\eps^2+|\alpha|^2}}$ and 
$\Vert C(\eps,\alpha)+i\eps A+\alpha B\Vert\leq c\,(\eps^2+|\alpha|^2)$ we note that the 
operator on the left hand side is a small perturbation of the identity. Therefore, this 
equation is uniquely solvable for all $(\eps,\alpha)\in M\cap B_\delta$ for sufficiently 
small $\delta>0$, and the solution $v$ depends continuously on $(\eps,\alpha)\in 
M\cap B_\delta$, and $\Vert v\Vert_H\leq c\Vert r\Vert_{C^1(\overline{B_\delta},H)}$ for all 
$(\eps,\alpha)\in M\cap B_\delta$ because $\Vert\tilde{u}(\eps,\alpha)\Vert_H\leq
\frac{c}{\sqrt{\eps^2+|\alpha|^2}}$ and $\Vert g(0,0)-g(\eps,\alpha)\Vert_H\leq 
c\sqrt{\eps^2+|\alpha|^2}$. This implies that $u^N:=\tilde{u}-v$ satisfies 
$C(\eps,\alpha)u^N(\eps,\alpha)=g(\eps,\alpha)$. Furthermore, in the case that all 
$\lambda_\ell$ have the same sign $\sigma$ the function $v(\eps,\cdot)$ is analytic in 
$\{\alpha\in\cmplx:|\alpha|<\delta,\ \sigma\Im\alpha<0\}$. If $\eps=0$ then $v(0,\alpha)$ 
is analytic and uniformly bounded in $\{\alpha\in\cmplx:|\alpha|<\delta,\ \alpha\not=0\}$. 
Finally, we have that 
\begin{eqnarray*}
u^R(\eps,\alpha) & = & T(\eps,\alpha)\,Q\,K(\eps,\alpha)u^N(\eps,\alpha)\ +\ 
T(\eps,\alpha)\,Qr(\eps,\alpha)\bigr) \\
& = & T(\eps,\alpha)\,Q\bigl(K(\eps,\alpha)-K(0,0)\bigr)u^N(\eps,\alpha)\ +\ 
T(\eps,\alpha)\,Qr(\eps,\alpha)\bigr)
\end{eqnarray*}
because $QK(0,0)u^N(\eps,\alpha)=Qu^N(\eps,\alpha)=0$. This shows continuity and boundedness 
of $u^R(\eps,\alpha)$ and ends the proof. \qed
\biglf
Later we only need the following conclusions from the previous theorem.
\begin{corollary} \label{c-abstract}
Let the assumptions of part (b) of the previous theorem hold. Then the following holds.
\begin{itemize}
\item[(a)] $I-K(\eps,\alpha)$ is an isomorphism from $H$ onto itself for all 
$(\eps,\alpha)\in M\cap B_\delta=\bigl\{(\eps,\alpha)\in[0,\delta)\times\cmplx:
|\alpha|<\delta,\ \sigma\Im\alpha\leq 0,\ \eps+|\alpha|>0\bigr\}$, and the mapping
$(\eps,\alpha)\mapsto\bigl(I-K(\eps,\alpha)\bigr)^{-1}$ is continuous from 
$M\cap B_\delta$ into $\cL(H,H)$.
\item[(b)] For fixed $\eps\in(0,\delta)$ and $r\in H$ the mapping $\alpha\mapsto
\bigl(I-K(\eps,\alpha)\bigr)^{-1}r$ is analytic on $\{\alpha\in\cmplx:|\alpha|<\delta,\ 
\sigma\Im\alpha<0\}$.
\end{itemize}
\end{corollary}
\biglf
We want to apply Theorem~\ref{t-absorp-abstract} and Corollary~\ref{c-abstract} to the 
equation~(\ref{eq-operator}) in a neighborhood $V\subset\cmplx$ of some propagative 
wavenumber $\hat{\alpha}_j$ for some $j\in J$ where $V$ is chosen such that it does not 
contain a cut-off value. (This is possible by Assumption~\ref{assump1} and 
Lemma~\ref{l-equiv}.) Therefore, we set $K(\eps,\alpha):=K_{\eps,\hat{\alpha}_j+\alpha}$ 
and $r(\eps,\alpha):=R_{\hat{\alpha}_j+\alpha}f$ and have to check the assumptions 
of the previous theorem. The operator $K$ and the right hand side $r$ are continuously 
differentiable with respect to $(\eps,\alpha)$ in a neighborhood of $(0,0)$. It has been 
shown in \cite{Kirsc2017} that the eigenvalue problem (\ref{ewp}) is equivalent to the 
eigenvalue problem (\ref{evp:a}), (\ref{evp:b}) and that all the other assumptions of 
the theorem are satisfied under the following additional assumption.
\begin{assumption} \label{assump2}
$\lambda_{\ell,j}\not=0$ for all $\ell=1,\ldots,m_j$ and $j\in J$.
\end{assumption}
\medlf
In \cite{Kirsc2017} we have shown that the application of Theorem~\ref{t-absorp-abstract} 
allows us to take the limit $\eps\to 0$. The following result has been shown.
\begin{theorem} \label{t-LAP1}
The solution $u_\eps$ converges to some solution $u$ of $\Delta u+k^2nu=-f$ in $H^1(K)$ for 
every $K$ of the form $K=(-R,R)\times(0,H)$. Furthermore, $u=0$ for $x_2=0$, and $u$ satisfies 
the radiation condition of Definition~\ref{d-RC} below for $q=0$ with
$$ \psi^\pm(x_1)\ =\ \frac{1}{2}\left[1\ \pm \frac{2}{\pi}\int\limits_0^{x_1/2}
\frac{\sin t}{t}\,dt\right],\quad x_1\in\real\,. $$
The coefficients $a_{\ell,j}$ in \eqref{decomp} are given by
\begin{equation} \label{eq:coeff_a}
a_{\ell,j}\ =\ \frac{2\pi i}{|\lambda_{\ell,j}|}\int\limits_Qf(x)\,
\overline{\hat{\phi}_{\ell,j}(x)}\,dx\,.
\end{equation}
\end{theorem}
\begin{definition} \label{d-RC}
Let $\psi^\pm\in C^\infty(\real)$ be any pair of functions with $\psi^\pm(x_1)=1+\cO(1/|x_1|)$
as $\pm x_1\to\infty$ and $\psi^\pm(x_1)=\cO(1/|x_1|)$ as $\pm x_1\to -\infty$ .
\newline
A solution $u$ of $\Delta u+k^2(n+q)u=-f$ in $\real^2_+$ satisfies the open waveguide 
radiation condition if $u$ has a decomposition into $u=u_{rad}+u_{prop}$ where:
\begin{itemize}
\item[(i)] The propagating part has the form
\begin{equation} \label{decomp}
u_{prop}(x)\ =\ \psi^+(x_1)\sum_{j\in J}\sum_{\ell:\lambda_{\ell,j}>0}a_{\ell,j}\,
\hat{\phi}_{\ell,j}(x)\ +\ \psi^-(x_1)\sum_{j\in J}\sum_{\ell:\lambda_{\ell,j}<0}a_{\ell,j}\,
\hat{\phi}_{\ell,j}(x)
\end{equation}
for some $a_{\ell,j}\in\cmplx$.
\item[(ii)] The radiating part satisfies $u_{rad}\in H^1(W_H)$ for all $H>h_0$ and its Fourier
transform $(\cF u_{rad})(\omega,x_2)$ with respect to $x_1$ satisfies the generalized 
Sommerfeld radiation condition
\begin{equation} \label{Sommerfeld}
\int\limits_{-\infty}^\infty\left|\partial_2(\cF u_{rad})(\omega,x_2)-i\sqrt{k^2-\omega^2}\,
(\cF u_{rad})(\omega,x_2)\right|^2d\omega\ \longrightarrow\ 0\,,\quad x_2\to\infty\,.
\end{equation}
Here we define the Fourier transform as
$$ (\cF\phi)(\omega)\ =\ \frac{1}{\sqrt{2\pi}}\int\limits_{-\infty}^\infty\phi(s)\,
e^{-is\omega}\,ds\,,\quad\omega\in\real\,. $$
\end{itemize}
\end{definition}
We note that we can replace $\psi^\pm$ by any pair of functions $\tilde{\psi}^\pm$ with 
$\tilde{\psi}^\pm(x_1)=1$ for $\pm x_1\geq T+1$ (for some $T>2\pi$) and 
$\tilde{\psi}^\pm(x_1)=0$ for $\pm x_1\leq T$. This is because the differences
$(\psi^\pm-\tilde{\psi}^\pm)\hat{\phi}_{\ell,j}$ are in $H^1(W_H)$ for all $H$ and decay 
exponentially as $x_2\to\infty$ and thus can be subsumed into the radiating part.
This representation give rise to the open waveguide radiation condition.
\biglf
In this paper we do not repeat the proof but use Corollary~\ref{c-abstract} in a different 
way. For the remaining part of the paper we make the following assumption.
\begin{assumption} \label{assump3}
For every $j\in J$ all of the eigenvalues $\lambda_{1,j},\ldots,\lambda_{m_j,j}$ have the 
same sign $\sigma_j\in\{+,-\}$.
\end{assumption}
Then we group the propagative wavenumber $\hat{\alpha}_j$ into right- and left going by 
defining $J^\pm:=\{j\in J:\sigma_j=\pm\}$. Furthermore, we define the following sets
\begin{eqnarray*}
\cA_j & := & \{\alpha\in\cmplx:|\alpha-\hat{\alpha}_j|<\delta,\ \pm\Im\alpha<0\}\quad
\mbox{ for }j\in J^\pm\quad\mbox{(open half discs)}\,, \\
\cA & := & I\cup\bigcup_{j\in J}\overline{\cA_j}\quad\mbox{where} \quad
I:=[-1/2,1/2]\setminus\bigcup_{j\in J}(\hat{\alpha}_j-\delta,\hat{\alpha}_j+\delta)\,, \\
\cM & := & \{(\eps,\alpha)\in[0,\delta]\times\cA:\eps+|\alpha-\hat{\alpha}_j|>0
\mbox{ for all }j\in J\}\,.
\end{eqnarray*}
\begin{figure}[h]
\begin{tikzpicture}[scale=1,thick]

\draw[thin,->] (0,-1) -- (0,1);
\draw[thin,->] (-6.5,0) -- (7,0);
\draw (7,-0.2) node{\tiny $\Re\alpha$};
\draw (0.2,1.2) node{\tiny $\Im\alpha$};
\filldraw[red] (2.3,0) arc (0:180:.3);
\filldraw[red] (4.7,0) arc (180:360:.3);
\filldraw[red] (-4.7,0) arc (0:180:.3);
\filldraw[red] (-2.3,0) arc (180:360:.3);
\draw[thin,-] (-6,-.05) -- (-6,.05);
\draw[thin,-] (6,-.05) -- (6,.05);
\draw (-6,-0.3) node{\tiny $-1/2$};
\draw (6,-0.3) node{\tiny $1/2$};
\draw[red,-] (-6,0) -- (6,0);
\filldraw[blue] (-5,0) circle (.05);
\filldraw[green] (-2,0) circle (.05);
\filldraw[blue] (2,0) circle (.05);
\filldraw[green] (5,0) circle (.05);
\end{tikzpicture}
\caption{The set $\cA$ (red), $\hat{\alpha}_j$ for $j\in J^+$ (green),  
$\hat{\alpha}_j$ for $j\in J^-$ (blue)}
\end{figure}
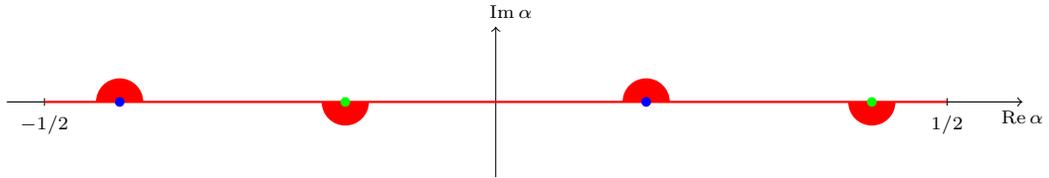
\biglf
Application of Corollary~\ref{c-abstract} to the equation (\ref{eq-operator}) yields the 
following result.
\begin{lemma} \label{l-main1}
Let Assumptions~\ref{assump1}, \ref{assump2}, and \ref{assump3} hold. Then there exists 
$\delta>0$ such that the operators $I-K_{\eps,\alpha}$, defined in (\ref{KR:a}), are 
isomorphisms from $H^1_{per,0}(Q)$ onto itself for all $(\eps,\alpha)\in\cM$, and 
$(\eps,\alpha)\mapsto\bigl(I-K_{\eps,\alpha}\bigr)^{-1}$ is continuous on $\cM$. If 
$\cK\subset\cM$ is compact then $\bigl\Vert\bigl(I-K_{\eps,\alpha}\bigr)^{-1}\bigr\Vert$ is 
uniformly bounded with respect to $(\eps,\alpha)\in\cK$. Furthermore, for every 
$\eps\in(0,\delta)$ the unique solution $\tilde{u}(\eps,\alpha)$ of (\ref{eq-operator}) 
depends analytically on $\alpha\in\bigcup_{j\in J}\cA_j$.
\smalf 
Formulated in terms of the variational problem (\ref{var-1}) we have the following. For 
every $f\in L^2(Q)$ and $(\eps,\alpha)\in\cM$ there exist a unique solution 
$\tilde{u}(\cdot,\eps,\alpha)\in H^1_{per,0}(Q)$ of (\ref{var-1}) which depends continuously 
on $(\eps,\alpha)\in\cM$; that is, $\tilde{u}\in C\bigl(\cM,H^1_{per,0}(Q)\bigr)$. For every 
compact set $\cK\subset\cM$ there exists $c>0$ with
$\Vert\tilde{u}(\cdot,\eps,\alpha)\Vert_{H^1(Q)}\leq c\Vert f\Vert_{L^2(Q)}$ for all
$(\eps,\alpha)\in\cK$. Finally, for every $\eps\in(0,\delta)$ the solution 
$\tilde{u}(\cdot,\eps,\alpha)$ depends analytically on $\alpha\in\bigcup_{j\in J}\cA_j$.
\end{lemma}
We recall that the original field $u_\eps$ is given by the inverse Floquet-Bloch transform; 
that is,
$$ u_\eps(x_1+2\pi p,x_2)\ =\ \int\limits_{-1/2}^{1/2}\tilde{u}(x,\eps,\alpha)\,
e^{i\alpha x_1}\,e^{2\pi p\alpha i}\,d\alpha\,,\quad x\in Q\,,\ p\in\ganz\,. $$
Therefore, if $\tilde{u}(x,\eps,\alpha)$ is analytic in, say, $\{\alpha\in\cmplx:
|\alpha-\hat{\alpha}_j|<\delta,\ \Im\alpha<0\}$ then we can modify the path of integration in 
this integral from the segment $(\hat{\alpha}_j-\delta/2,\hat{\alpha}_j+\delta/2)$ to 
$\Gamma_j^-:=\{\alpha\in\cmplx:|\alpha-\hat{\alpha}_j|=\delta/2,\ \Im\alpha<0\}$. With the 
analogous change for $j\in J^-$ we define the global path 
$$ \Gamma\ :=\ I \cup \bigcup_{j\in J^+}\Gamma_j^- \cup \bigcup_{j\in J^-}\Gamma_j^+
\quad\mbox{where now}\quad I:=[-1/2,1/2]\setminus\bigcup_{j\in J}
(\hat{\alpha}_j-\delta/2,\hat{\alpha}_j+\delta/2) $$
which connects $-1/2$ with $+1/2$. Application of the previous lemma yields the following 
form of the Limiting Absorption Principle.
\begin{theorem} \label{LAP1}
Let Assumptions~\ref{assump1}, \ref{assump2}, and \ref{assump3} hold. Then there exists 
$\delta>0$ such that for all $f\in L^2(Q)$ and $\eps\in(0,\delta)$ there exist a unique 
solution $u_\eps\in H^1_0(\real^2_+)$ of (\ref{source-1}) for $q=0$ and $k$ replaced by 
$k+i\eps$, and $u_\eps$ has the form
$$ u_\eps(x_1+2\pi p,x_2)\ =\ \int\limits_\Gamma\tilde{u}(x,\eps,\alpha)\,e^{i\alpha x_1}\,
e^{2\pi p\alpha i}\,d\alpha\,,\quad x\in Q\,,\ p\in\ganz\,, $$
where $\tilde{u}(\cdot,\eps,\alpha)\in H^1_{per,0}(Q)$ satisfies (\ref{var-1}) for all 
$\alpha\in\Gamma$. Furthermore, $u_\eps$ converges to some solution 
$u\in H^1_{loc,0}(W_{h_0})$ of (\ref{source-1}) for $q=0$. The convergence is in $H^1(K)$ 
for every $K=(-R,R)\times(0,h_0)$. The function $u$ has the representation as
$$ u(x_1+2\pi p,x_2)\ =\ \int\limits_\Gamma\tilde{u}(x,0,\alpha)\,e^{i\alpha x_1}\,
e^{2\pi p\alpha i}\,d\alpha\,,\quad x\in Q\,,\ p\in\ganz\,, $$
where $\tilde{u}(\cdot,0,\alpha)\in H^1_{per,0}(Q)$ is the solution of (\ref{var-1}) for
$\alpha\in\Gamma$ and $\eps=0$.
\end{theorem}
\textbf{Proof:} This is an immediate consequence of the fact that 
$\tilde{u}(\cdot,\eps,\alpha)$ depends continuously on $(\eps,\alpha)\in[0,\delta)\times
\Gamma\subset\cM$. \qed
\biglf
\textbf{Remark:} We note that uniqueness holds for (\ref{var-2}) and (\ref{eq-operator}) for 
$\eps=0$ because $I-K_{0,\alpha}$ is an isomorphism for all $\alpha\in\Gamma$.
\biglf
We have to furthermore modify the contour $\Gamma$ in neighborhoods of the cut-off values
because the square root function in the definition (\ref{KR:a}) of the operator 
$K_{0,\alpha}$ does not depend analytically on $\alpha$. Let $k=\hat{\ell}+\kappa$ with 
$\hat{\ell}\in\ganz_{\geq 0}$ and $\kappa\in(-1/2,1/2]$. We make the additional assumption.

\begin{assumption}\label{assump4}
Assume that the wavenumber $k>0$ satisfies $k\notin\frac{1}{2}\nat$.
\end{assumption}
 
 With Assumption \ref{assump4}, $0<|\kappa|<1/2$ thus cut-off values are given by 
$\pm\kappa$. We show how to modify the path $\Gamma$ of integration in the neighborhoods 
of $\pm\kappa$. Recall that $\sqrt{k^2-(\ell+\alpha)^2}=\sqrt{\hat{\ell}+\kappa-\ell-\alpha}
\sqrt{\hat{\ell}+\kappa+\ell+\alpha}$ We decompose
\begin{eqnarray*}
& & \sum_{\ell\in\ganz}\sqrt{k^2-(\ell+\alpha)^2}\,v_\ell(h_0)\,\overline{\psi_\ell(h_0)}\ =\
\sum_{\ell\not=\pm\hat{\ell}}\sqrt{k^2-(\ell+\alpha)^2}\,v_\ell(h_0)\,
\overline{\psi_\ell(h_0)} \\
& & +\ \sqrt{\kappa-\alpha}\,\sqrt{2\hat{\ell}+\kappa+\alpha}\,v_{\hat{\ell}}(h_0)\,
\overline{\psi_{\hat{\ell}}(h_0)}\ +\
\sqrt{\kappa+\alpha}\,\sqrt{2\hat{\ell}+\kappa-\alpha}\,v_{-\hat{\ell}}(h_0)\,
\overline{\psi_{-\hat{\ell}}(h_0)}\,.
\end{eqnarray*}
Since $\sqrt{\kappa-\alpha}\,\sqrt{2\hat{\ell}+\kappa+\alpha}$ is holomorphic in 
$\{\alpha\in\cmplx:|\alpha-\kappa|<\delta,\ \Im\alpha<0\}$ for sufficiently small $\delta>0$ 
we modify $\Gamma$ and replace $(\kappa-\delta/2,\kappa+\delta/2)$  by the half circle 
$\{\alpha\in\cmplx:|\alpha-\kappa|=\delta/2,\ \Im\alpha<0\}$ in the definition of $\Gamma$. 
Analogously, we replace $(-\kappa-\delta/2,-\kappa+\delta/2)$  by the half circle 
$\{\alpha\in\cmplx:|\alpha+\kappa|=\delta/2,\ \Im\alpha>0\}$. In the same way we replace the 
set $\cA$ by $\cA\cup\{\alpha\in\cmplx:|\alpha-\kappa|<\delta,\ \Im\alpha<0\}\cup
\{\alpha\in\cmplx:|\alpha+\kappa|<\delta,\ \Im\alpha>0\}$.
 
\section{The Perturbed Case}

We recall the definition of the path $\Gamma$, observe that $\Gamma\subset\cA$, 
and define the operator $P$ from $C\bigl(\cA,H^1_{per,0}(Q)\bigr)$ into $L^2(Q)$ by 
$$ (Pg)(x)\ :=\ \int\limits_\Gamma g(x,\beta)\,e^{i\beta x_1}\,d\beta\,,\quad x\in Q\,,
\quad g\in C\bigl(\cA,H^1_{per,0}(Q)\bigr)\,. $$
Then it is easily seen that $P$ is compact. Indeed, let 
$g_j\in C\bigl(\cA,H^1_{per,0}(Q)\bigr)$ be bounded; that is, 
$\sup_{\beta\in\cA}\Vert g_j(\cdot,\beta)\Vert_{H^1(Q)}\leq c$ for all $j\in\nat$. Then 
$\Vert Pg_j\Vert_{H^1(Q)}\leq c_1\int_\Gamma\Vert g_j(\cdot,\beta)\Vert_{H^1(Q)}d\beta
\leq c_2$ for all $j$. The compactness of $H^1_{per,0}(Q)$ in $L^2(Q)$ yields compactness 
of $P$.
\medlf
We consider now the local perturbation of the periodic case; that is, we look at 
(\ref{source-1}) for arbitrary $q\in L^2(Q)$. We rewrite the equation for $k+i\eps$ instead 
of $k$ as
$$ \Delta u_\eps\ +\ (k+i\eps)^2\,n\,u_\eps\ =\ -f\ -\ (k+i\eps)^2\,q\,u_\eps\quad\mbox{in }
\real^2_+\,. $$
The (periodic) Floquet-Bloch transformed equation has the form 
\begin{eqnarray}
& & \Delta\tilde{u}_\eps(\cdot,\alpha)\ +\ 2i\alpha\,\partial_1\tilde{u}_\eps(\cdot,\alpha)\ +\
[(k+i\eps)^2\,n-\alpha^2]\,\tilde{u}_\eps(\cdot,\alpha) \nonumber \\ 
& = & -f\,e^{-i\alpha x_1}\ -\ (k+i\eps)^2\,q\,u_\eps\,e^{-i\alpha x_1} \nonumber \\ 
& = & -f\,e^{-i\alpha x_1}\ -\ (k+i\eps)^2\,q\,e^{-i\alpha x_1}\int\limits_{-1/2}^{1/2}
\tilde{u}_\eps(\cdot,\beta)\,e^{i\beta x_1}\,d\beta\quad\mbox{in }Q^\infty \label{eq-pert-1}
\end{eqnarray}
for $\alpha\in[-1/2,1/2]$ and in variational form (compare with (\ref{var-2}))
\begin{eqnarray}
& & \int\limits_Q\biggl[\nabla\tilde{u}_\eps(\cdot,\alpha)\cdot\nabla\overline{\psi}-2i\alpha\,
\partial_1\tilde{u}_\eps(\cdot,\alpha)\,\overline{\psi}-[(k+i\eps)^2\,n-\alpha^2]\,
\tilde{u}_\eps(\cdot,\alpha)\,\overline{\psi}\biggr]\,dx\ - \label{var-3} \\
& &  -\ i\sum_{\ell\in\ganz}\sqrt{(k+i\eps)^2-(\ell+\alpha)^2}\,
\tilde{u}_{\eps,\ell}(h_0,\alpha)\,\overline{\psi_\ell(h_0)} \nonumber \\ 
& = & \int\limits_Q[f\,e^{-i\alpha x_1}+(k+i\eps)^2\,q\,u_\eps]\,\overline{\psi}\,dx 
\nonumber \\
& = & \int\limits_Qf\,e^{-i\alpha x_1}\,\overline{\psi}\,dx\ +\ (k+i\eps)^2
\int\limits_Q \,q\,e^{-i\alpha x_1}\,\overline{\psi}\int\limits_{-1/2}^{1/2}
\tilde{u}_\eps(\cdot,\beta)\,e^{i\beta x_1}\,d\beta\,dx \nonumber
\end{eqnarray}
for all $\psi\in H^1_{per,0}(Q)$. As in (\ref{var-2}) we write this as
$$ (I-K_{\eps,\alpha})\tilde{u}_\eps(\cdot,\alpha)\ =\ R_\alpha f\ +\ (k+i\eps)^2\,
R_\alpha\left(q\int\limits_{-1/2}^{1/2}\tilde{u}_\eps(\cdot,\beta)\,e^{i\beta x_1}d\beta
\right)\,. $$
Now we observe that the operators $K_{\eps,\alpha}$ depend smoothly on $\alpha\in\cA$. 
Considering the right hand side as a source and modifying the path of integration is the 
motivation to study (for fixed $\eps>0$) the equation
\begin{equation} \label{eq-pert-2}
(I-K_{\eps,\alpha})v_\eps(\cdot,\alpha)\ =\ R_\alpha f\ +\ (k+i\eps)^2\,
R_\alpha\bigl(q\,Pv_\eps\bigr)\quad\mbox{in }H^1_{per,0}(Q)
\end{equation}
for $\alpha\in\Gamma$ where we recall the definitions of $K_{\eps,\alpha}:H^1_{per,0}(Q)\to 
H^1_{per,0}(Q)$ and $R_\alpha:L^2(Q)\to H^1_{per,0}(Q)$ from (\ref{KR:a}), (\ref{KR:b}). 
The corresponding variational formulation is (\ref{var-3}) with $\int_{-1/2}^{1/2}
v_\eps(\cdot,\beta)\,e^{i\beta x_1}\,d\beta$ replaced by $Pv_\eps$; that is,
\begin{eqnarray}
& & \int\limits_Q\biggl[\nabla v_\eps(\cdot,\alpha)\cdot\nabla\overline{\psi}-2i\alpha\,
\partial_1v_\eps(\cdot,\alpha)\,\overline{\psi}-[(k+i\eps)^2\,n-\alpha^2]\,
v_\eps(\cdot,\alpha)\,\overline{\psi}\biggr]\,dx\ - \label{var-4} \\
& &  -\ i\sum_{\ell\in\ganz}\sqrt{(k+i\eps)^2-(\ell+\alpha)^2}\,
v_{\eps,\ell}(h_0,\alpha)\,\overline{\psi_\ell(h_0)} \nonumber \\ 
& = & \int\limits_Qf\,e^{-i\alpha x_1}\,\overline{\psi}\,dx\ +\ (k+i\eps)^2
\int\limits_Q \,q\,e^{-i\alpha x_1}\,\overline{\psi}\int\limits_\Gamma
v_\eps(\cdot,\beta)\,e^{i\beta x_1}\,d\beta\,dx \nonumber
\end{eqnarray}
for all $\psi\in H^1_{per,0}(Q)$ and $\alpha\in\Gamma$.

\begin{lemma} \label{l-main2} 
Let Assumptions~\ref{assump1}, \ref{assump2}, \ref{assump3},  and \ref{assump4} hold. Then there exists 
$\delta>0$ such that for all $f\in L^2(Q)$ and $\eps\in(0,\delta)$ there exist a unique 
solution $v_\eps\in C\bigl(\Gamma,H^1_{per,0}(Q)\bigr)$ of (\ref{eq-pert-2}) and 
(\ref{var-4}). The function $(\eps,\alpha)\mapsto v_\eps(\cdot,\alpha)$ can be extended to 
$v\in C\bigl(\cM,H^1_{per,0}(Q)\bigr)$ and $v(\cdot,\eps,\alpha)$ depends analytically on 
$\alpha\in\bigcup_{j\in J}\cA_j$ for every $\eps\in(0,\delta)$.
\newline
Furthermore, $v(\cdot,\eps,\alpha)\in H^1_{per,0}(Q)$ is the unique solution of 
(\ref{eq-pert-1}) for all $(\eps,\alpha)\in\cM$.
\end{lemma}
\textbf{Proof:} First we introduce the operator $T:L^2(Q)\to 
C\bigl(\Gamma,H^1_{per,0}(Q)\bigr)$ as $Tg=w_g$ where 
$w_g\in C\bigl(\Gamma,H^1_{per,0}(Q)\bigr)$ is the unique solution of 
$$ (I-K_{\eps,\alpha})w_g(\cdot,\alpha)\ =\ R_\alpha g\quad\mbox{in }H^1_{per,0}(Q)\quad
\mbox{for all }\alpha\in\Gamma\,. $$
The existence and boundedness of $T$ is assured by Lemma~\ref{l-main1}. Then we can write 
(\ref{eq-pert-2}) as the fixpoint equation 
$$ v\ =\ T\bigl(f+(k+i\eps)^2q\,Pv\bigr)\ =\ Tf\ +\ (k+i\eps)^2\,T(q\,Pv)\quad\mbox{in }
C\bigl(\Gamma,H^1_{per,0}(Q)\bigr) $$ 
which is of Fredholm type because of the compactness of $P$. Therefore, it suffices 
to prove uniqueness. Let $v\in C\bigl(\Gamma,H^1_{per,0}(Q)\bigr)$ be a solution corresponding 
to $f=0$. Define $w(\cdot,\alpha)\in H^1_{per,0}(Q)$ as the unique solution of 
$$ (I-K_{\eps,\alpha})w(\cdot,\alpha)\ =\ (k+i\eps)^2\,R_\alpha\bigl(q\,Pv\bigr)\quad
\mbox{in }H^1_{per,0}(Q)\quad\mbox{for all }\alpha\in\cA\,. $$
Then $w\in C\bigl(\cA,H^1_{per,0}(Q)\bigr)$, again by Lemma~\ref{l-main1}. Furthermore,
$(I-K_{\eps,\alpha})\bigl(v(\cdot,\alpha)-w(\cdot,\alpha)\bigr)=0$ for all $\alpha\in\Gamma$
which implies that $v(\cdot,\alpha)=w(\cdot,\alpha)\bigr)$ for all $\alpha\in\Gamma$. 
Therefore, $Pv=Pw$ which shows that $w$ satisfies
\begin{equation} \label{aux1}
(I-K_{\eps,\alpha})w(\cdot,\alpha)\ =\ (k+i\eps)^2\,R_\alpha\bigl(q\,Pw\bigr)\quad
\mbox{in }H^1_{per,0}(Q)\quad\mbox{for all }\alpha\in\cA\,. 
\end{equation}
Since $w$ depends analytically on $\alpha\in\bigcup_{j\in J}\cA_j$ we can change the path 
of integration and have 
$$ (Pw)(x)\ :=\ \int\limits_{-1/2}^{1/2}w(x,\beta)\,e^{i\beta x_1}\,d\beta\,,\quad x\in Q\,. $$
Therefore, (\ref{aux1}) implies that $w(\cdot,\alpha)$ satisfies (\ref{var-3}) for $f=0$ and 
all $\alpha\in[-1/2,1/2]$. Standard arguments (setting $\psi=w$ and taking the imaginary part)
yields that $w(\cdot,\alpha)$ vanishes for all $\alpha\in[-1/2,1/2]$. Therefore, $Pw=0$ 
and $w(\cdot,\alpha)$ satisfies $(I-K_{\eps,\alpha})w(\cdot,\alpha)=0$ for all $\alpha\in\cA$ 
which implies again that $w(\cdot,\alpha)$ vanishes in $Q$ for all $\alpha\in\cA$.
\newline 
So far, we kept $\eps>0$ fixed. However, the extension $w=w(\cdot,\eps,\alpha)$ of 
$v=v_\eps(\cdot,\alpha)$ is in $C\bigl(\cM,H^1_{per,0}(Q)\bigr)$ which ends the proof. \qed
\biglf
Since $(\eps,\alpha)\mapsto v(\cdot,\eps,\alpha)$ is continuous on $[0,\delta)\times\Gamma$ we 
can let $\eps$ tend to zero in (\ref{var-4}). We have therefore shown the Limiting Absorption 
Principle in the following form.
\begin{theorem} \label{LAP2}
Let Assumptions~\ref{assump1}, \ref{assump2},  \ref{assump3} and \ref{assump4} hold. Then there exists 
$\delta>0$ such that for all $f\in L^2(Q)$ and $\eps\in(0,\delta)$ there exist a unique 
solution $u_\eps\in H^1_0(\real^2_+)$ of (\ref{source-1}) for $k$ replaced by $k+i\eps$, and 
$u_\eps$ has the form
$$ u_\eps(x_1,x_2)\ =\ \int\limits_\Gamma v_\eps(x,\alpha)\,e^{i\alpha x_1}\,d\alpha\,,
\quad x\in W_{h_0}\,, $$
where $v_\eps(\cdot,\alpha)\in H^1_{per,0}(Q)$ satisfies (\ref{var-4}) for all 
$\alpha\in\Gamma$. Furthermore, $u_\eps$ converges to some solution 
$u\in H^1_{loc,0}(W_{h_0})$ of (\ref{source-1}). The convergence is in $H^1(K)$ for every
$K=(-R,R)\times(0,h_0)$. The function $u$ has the representation as
\begin{equation} \label{u-def}
u(x_1,x_2)\ =\ \int\limits_\Gamma v_0(x,\alpha)\,e^{i\alpha x_1}\,\,d\alpha\,,\quad 
x\in W_{h_0}\,, 
\end{equation}
where $v_0(\cdot,\alpha)\in H^1_{per,0}(Q)$ is the solution of (\ref{var-4}) for
$\alpha\in\Gamma$ and $\eps=0$; that is,
\begin{equation} \label{eq-operator-0}
(I-K_{0,\alpha})v_0(\alpha)\ =\ R_\alpha f\ +\ k^2\,R_\alpha(q\,Pv_0)\quad\mbox{in }
H^1_{per,0}(Q)\,,\ \alpha\in\Gamma\,,
\end{equation}
or, in variational form,
\begin{eqnarray}
& & \int\limits_Q\biggl[\nabla v_0(\cdot,\alpha)\cdot\nabla\overline{\psi}-2i\alpha\,
\partial_1v_0(\cdot,\alpha)\,\overline{\psi}-[k^2\,n-\alpha^2]\,
v_0(\cdot,\alpha)\,\overline{\psi}\biggr]\,dx\ - \label{var-5} \\
& &  -\ i\sum_{\ell\in\ganz}\sqrt{k^2-(\ell+\alpha)^2}\,
v_{0,\ell}(h_0,\alpha)\,\overline{\psi_\ell(h_0)} \nonumber \\ 
& = & \int\limits_Qf\,e^{-i\alpha x_1}\,\overline{\psi}\,dx\ +\ k^2
\int\limits_Q \,q\,e^{-i\alpha x_1}\,\overline{\psi}\int\limits_\Gamma
v_0(\cdot,\beta)\,e^{i\beta x_1}\,d\beta\,dx \nonumber
\end{eqnarray}
for all $\psi\in H^1_{per,0}(Q)$ and $\alpha\in\Gamma$. $v_0$ is continuous with 
respect to $\alpha$ and thus $v_0\in C\bigl(\Gamma,H^1_{per,0}(Q)\bigr)$.
\end{theorem}
\begin{remark}
(\ref{var-5}) is the variational formulation of the following coupled system of boundary 
value problems (formulated only for $\eps=0$) for $v_0(\cdot,\alpha)\in H^1_{per,0}(Q)$:
\begin{eqngroup}\begin{eqnarray}
& & \Delta v_0(\cdot,\alpha)\ +\ 2i\alpha\,\partial_1v_0(\cdot,\alpha)\ +\
[k^2\,n-\alpha^2]\,v_0(\cdot,\alpha) \nonumber \\ 
& = & -f\,e^{-i\alpha x_1}\ -\ k^2\,q\,e^{-i\alpha x_1}\int\limits_\Gamma
v_0(\cdot,\beta)\,e^{i\beta x_1}\,d\beta\quad\mbox{in }Q\,, \label{eq_de:a} \\
v_0(\cdot,\alpha) & = & 0\quad\mbox{for }x_2=0\,, \label{eq_de:b} \\
\partial_2v_{0,\ell}(h_0,\alpha) & = & i\sqrt{k^2-(\ell+\alpha)^2}\,v_{0,\ell}(h_0,\alpha)\,, 
\label{eq_de:c}
\end{eqnarray} \end{eqngroup}
for all $\alpha\in\Gamma$.
\end{remark}
\begin{theorem} \label{t-RC}
Let $v_0\in C\bigl(\Gamma,H^1_{per,0}(Q)\bigr)$ be a solution of the system (\ref{var-5}). Then 
the corresponding field $u$, defined in (\ref{u-def}), has an extension into $\real^2_+$ which 
satisfies (\ref{source-1}) and the radiation condition.
\end{theorem}
\textbf{Proof:} From (\ref{eq-operator-0}) we conclude that 
$(I-K_{0,\alpha})v_0(\alpha)=g_\alpha$ in $H^1_{per,0}(Q)$ for all $\alpha\in\Gamma$ with 
$g_\alpha:=R_\alpha f+k^2R_\alpha\bigl(qPv_0(\alpha)\bigr)$. By Lemma~\ref{l-main1} 
$\alpha\mapsto v_0(\alpha)$ is continuous in $\cA$ and can be extended analytically into 
$\cup_{j\in J}\cA_j$. We consider the parts $\Gamma_j^\pm$ of the integration curve in 
(\ref{u-def}). Let, for example, $j\in J^+$. Then $v_0(x,\alpha)e^{i\alpha x_1}$ has the 
representation (cf Theorem~\ref{t-absorp-abstract})
\begin{equation} \label{repres}
v_0(x,\alpha)\,e^{i\alpha x_1}\ =\ v^b(x,0,\alpha)\,e^{i\alpha x_1}\ +\ 
\frac{1}{\alpha-\hat{\alpha}_j}\sum_{\ell=1}^{m_j}\frac{r_{\ell,j}}
{\lambda_{\ell,j}}\,\hat{\phi}_{\ell,j}(x)\,e^{i(\alpha-\hat{\alpha}_j)x_1}\,,\quad x\in Q\,,
\end{equation}
for $0<|\alpha-\hat{\alpha}_j|<\delta$. The part $v^b$ is analytic in $\{\alpha\in\cmplx:
|\alpha-\hat{\alpha}_j|<\delta\}$ and $v_0$ analytic in $\{\alpha\in\cmplx:
0<|\alpha-\hat{\alpha}_j|<\delta\}$. Now we extend $v_0(\cdot,\alpha)$ into $Q^\infty\setminus 
Q$ by
$$ v_0(x,\alpha)\ :=\ \sum_{\ell\in\ganz}v_{0,\ell}\,
e^{i\sqrt{k^2-(\ell+\alpha)^2}(x_2-h_0)}\,e^{i\ell x_1}\,,\quad x_2>h_0\,, $$
for all $\alpha\in I\cup\{\alpha\in\cmplx:0<|\alpha-\hat{\alpha}_j|<\delta/2\}$. We note that 
for $|\ell|\geq k+1$ we have that $\Re[k^2-(\ell+\alpha)^2]=k^2-(\ell+\Re\alpha)^2+
(\Im\alpha)^2<0$ for sufficiently small $\delta>0$. Therefore, $\Im\sqrt{k^2-(\ell+\alpha)^2}
>0$, and the series $\sum_{|\ell|\geq k+1}v_{0,\ell}\,
e^{i\sqrt{k^2-(\ell+\alpha)^2}(x_2-h_0)}\,e^{i\ell x_1}$ decays exponentially as 
$x_2\to\infty$. The finite sum 
$\sum_{|\ell|<k+1}v_{0,\ell}\,e^{i\sqrt{k^2-(\ell+\alpha)^2}(x_2-h_0)}\,e^{i\ell x_1}$ can 
grow exponentially.
\newline 
Since the functions $\hat{\phi}_{\ell,j}$ are already solutions of the Helmholtz equation in 
all of $\real^2_+$ we can take (\ref{repres}) as the definition of $v^b(x,0,\alpha)$ for 
$x_2>h_0$ and $0<|\alpha-\hat{\alpha}_j|<\delta$. Then $v^b$ is still analytic in 
$\{\alpha\in\cmplx:0<|\alpha-\hat{\alpha}_j|<\delta\}$. Therefore,
\begin{eqnarray*}
\int\limits_{C_j^-}v_0(x,\alpha)\,e^{i\alpha x_1}\,d\alpha & = &
\int\limits_{C_j^-}u^b(x,0,\alpha)\,e^{i\alpha x_1}\,d\alpha\ +\ 
\int\limits_{C_j^-}\frac{1}{\alpha-\hat{\alpha}_j}\,e^{i(\alpha-\hat{\alpha}_j)x_1}\,d\alpha
\sum_{\ell=1}^{m_j}\frac{r_{\ell,j}}{\lambda_{\ell,j}}\,
\hat{\phi}_{\ell,j}(x) \\
& = &  \int\limits_{\hat{\alpha}_j-\delta/2}^{\hat{\alpha}_j+\delta/2}u^b(x,0,\alpha)
\,e^{i\alpha x_1}\,d\alpha\ +\ 
\int\limits_{C^-}\frac{1}{\alpha}\,e^{i\alpha x_1}\,d\alpha
\sum_{\ell=1}^{m_j}\frac{r_{\ell,j}}{\lambda_{\ell,j}}\,\hat{\phi}_{\ell,j}(x)
\end{eqnarray*}
for $x\in Q^\infty$ where $C^-:=\{\alpha\in\cmplx:|\alpha|=\delta/2,\ \Im\alpha<0\}$. We 
compute
\begin{eqnarray*}
\int\limits_{C^-}\frac{1}{\alpha}\,e^{i\alpha x_1}\,d\alpha & = & 
\int\limits_{C^-}\frac{d\alpha}{\alpha}\ +\ 
\int\limits_{C^-}\frac{\cos(\alpha x_1)-1}{\alpha}\,d\alpha\ +\ 
i\int\limits_{C^-}\frac{\sin(\alpha x_1)}{\alpha}\,d\alpha \\
& = & \pi\,i\ +\ \int\limits_{-\delta/2}^{\delta/2}\frac{\cos(\alpha x_1)-1}{\alpha}\,d\alpha\
+\ i\int\limits_{-\delta/2}^{\delta/2}\frac{\sin(\alpha x_1)}{\alpha}\,d\alpha
\end{eqnarray*}
because the integrands of the second and third integrals are holomorphic. The second integral 
vanishes because the integrand is an odd function. Substituting $t=\alpha x_1$ in the third 
integral yields
$$ \int\limits_{C^-}\frac{d\alpha}{\alpha}\ =\ \pi\,i\ +\ 2i\int\limits_0^{\delta x_1/2}
\frac{\sin t}{t}\,dt\,. $$
Integration for $j\in J^-$; that is, integration over $\Gamma_j^+$ is done analogously. 
Furthermore, integration over the half circles with centers $\pm\kappa$ are reduced to the 
integrals over the segments $(\pm\kappa-\delta/2,\pm\kappa+\delta/2)$. Therefore,
\begin{eqnarray*}
u(x) & = & \int\limits_\Gamma v_0(x,\alpha)\,e^{i\alpha x_1}d\alpha\ =\
\int\limits_I v_0(x,\alpha)\,e^{i\alpha x_1}d\alpha \\
& & +\ \sum_{j\in J^+}\int\limits_{\Gamma_j^-}v_0(x,\alpha)\,e^{i\alpha x_1}d\alpha\ +\
\sum_{j\in J^-}\int\limits_{\Gamma_j^+}v_0(x,\alpha)\,e^{i\alpha x_1}d\alpha \\
& = & \int\limits_I v_0(x,\alpha)\,e^{i\alpha x_1}d\alpha\ +\
\sum_{j\in J}\int\limits_{\hat{\alpha}_j-\delta/2}^{\hat{\alpha}_j+\delta/2}
v^b(x,0,\alpha)\,e^{i\alpha x_1}d\alpha \\
& & +\ 2\pi i\left[\frac{1}{2}+\frac{1}{\pi}\int\limits_0^{\delta x_1/2}
\frac{\sin t}{t}\,dt\right]\sum_{j\in J^+}\sum_{\ell=1}^{m_j}
\frac{r_{\ell,j}}{\lambda_{\ell,j}}\,\hat{\phi}_{\ell,j}(x) \\
& & -\ 2\pi i\left[\frac{1}{2}-\frac{1}{\pi}\int\limits_0^{\delta x_1/2}
\frac{\sin t}{t}\,dt\right]\sum_{j\in J^-}\sum_{\ell=1}^{m_j}
\frac{r_{\ell,j}}{\lambda_{\ell,j}}\,\hat{\phi}_{\ell,j}(x)\,,\quad x\in\real^2_+\,.
\end{eqnarray*}
With $\psi^\pm(x_1)=\frac{1}{2}\pm\frac{1}{\pi}\int_0^{\delta x_1/2}\frac{\sin t}{t}\,dt$ we
have shown the decomposition of $u$ into a radiating part $u_{rad}$ and the propagating part
$u_{prop}$ the the layer $W_{h_0}$.
\smalf 
We set
\begin{equation} \label{v_rad}
v_{rad}(\cdot,\alpha)\ :=\ \left\{\begin{array}{cl} v_0(\cdot,\alpha)\,, & \alpha\in I\,, \\
v^b(\cdot,0,\alpha)\,, & \alpha\in\bigcup_{j\in J}
(\hat{\alpha}_j-\delta/2,\hat{\alpha}_j+\delta/2)\,, \\
0\,, & \mbox{else}\,. \end{array}\right.
\end{equation}
Then $v_{rad}\in L^2\bigl((-1/2,1/2),H^1_{per,0}(Q^H)\bigr)$ for all $H>h_0$ where 
$Q^H=(0,2\pi)\times(0,H)$ and $u_{rad}(x)=\int_{-1/2}^{1/2}v_{rad}(x,\alpha)e^{i\alpha x_1}
d\alpha$. We conclude that $u_{rad}\in H^1(W_H)$ for all $H<h_0$ by the well known properties 
of the Floquet-Bloch transform. It remains to show the generalized angular spectral radiation 
condition (\ref{Sommerfeld}). First we note that, with $\omega=\ell+\alpha$ where 
$\ell\in\ganz$ and $\alpha\in(-1/2,1/2]$,
\begin{eqnarray*}
(\cF u_{rad})(\omega,x_2) & = & \frac{1}{\sqrt{2\pi}}\int\limits_{-\infty}^\infty 
u_{rad}(x_1,x_2)\,e^{-i\omega x_1}dx_1 \\
& = & \frac{1}{\sqrt{2\pi}}\sum_{m\in\ganz}\int\limits_0^{2\pi}u_{rad}(x_1+2\pi m,x_2)\,
e^{-i(\ell+\alpha)(x_1+2\pi m)}dx_1 \\
& = & \frac{1}{\sqrt{2\pi}}\int\limits_0^{2\pi}v_{rad}(x_1,x_2,\alpha)\,
e^{-i\ell x_1}dx_1\ =\ v_{rad,\ell}(x_2,\alpha)
\end{eqnarray*}
which are the Fourier coefficients of $x_1\mapsto v_{rad}(x,\alpha)$. Therefore,
\begin{eqnarray*}
& & \int\limits_{-\infty}^\infty\left|\partial_2(\cF u_{rad})(\omega,x_2)-
i\sqrt{k^2-\omega^2}\,(\cF u_{rad})(\omega,x_2)\right|^2\,d\omega \\
& = & \sum_{m\in\ganz}\int\limits_{-1/2}^{1/2}\left|\partial_2(\cF u_{rad})(\ell+\alpha,x_2)-
i\sqrt{k^2-(\ell+\alpha)^2}\,(\cF u_{rad})(\ell+\alpha,x_2)\right|^2\,d\alpha \\
& = & \sum_{m\in\ganz}\int\limits_{-1/2}^{1/2}\left|\partial_2v_{rad,\ell}(x_2,\alpha)-
i\sqrt{k^2-(\ell+\alpha)^2}\,v_{rad,\ell}(x_2,\alpha)\right|^2\,d\alpha\,.
\end{eqnarray*}
From the definition~(\ref{v_rad}) of $v_{rad}$ we observe that 
$\partial_2v_{rad,\ell}(x_2,\alpha)-i\sqrt{k^2-(\ell+\alpha)^2}\,v_{rad,\ell}(x_2,\alpha)=
\partial_2v_{0,\ell}(x_2,\alpha)-i\sqrt{k^2-(\ell+\alpha)^2}\,
v_{0,\ell}(x_2,\alpha)=0$ for $\alpha\in I$ and $x_2>h_0$. For $|\alpha-\hat{\alpha}_j|<
\delta/2$ with $j\in J^+$ we change the path of integration into $C_j^\mp$ and have from 
(\ref{repres}) that
\begin{eqnarray*}
& & \sum_{\ell\in\ganz}\int\limits_{C_j^\mp}\left|\partial_2v^b_\ell(x_2,0,\alpha)-
i\sqrt{k^2-(\ell+\alpha)^2}\,v^b_\ell(x_2,0,\alpha)\right|^2\,d\alpha \\
& \leq & c\sum_{\ell\in\ganz}\sum_{\ell^\prime=1}^{m_j}\left[
\left|\partial_2\hat{\phi}_{\ell^\prime,j,\ell}(x_2)\right|^2+
|\ell|^2\left|\hat{\phi}_{\ell^\prime,j,\ell}(x_2)\right|^2\right] \\
& \leq & c\sum_{\ell^\prime=1}^{m_j}\int\limits_0^{2\pi}\left[
\left|\partial_2\hat{\phi}_{\ell^\prime,j}(x_1,x_2)\right|^2+
\left|\partial_1\hat{\phi}_{\ell^\prime,j}(x_1,x_2)\right|^2\right]\,dx_1
\end{eqnarray*}
where $\hat{\phi}_{\ell^\prime,j,\ell}(x_2)=\frac{1}{\sqrt{2\pi}}\int_0^{2\pi}
\hat{\phi}_{\ell^\prime,j}(x_1,x_2)e^{-i(\ell+\hat{\alpha}_j)x_1}dx_1$ are the Fourier 
coefficients of \\ $\hat{\phi}_{\ell^\prime,j}(x_1,x_2)e^{-i\hat{\alpha}_jx_1}$. The 
exponential decay of $\hat{\phi}_{\ell^\prime,j}$ and its derivatives yields that this 
expression tends to zero as $x_2\to\infty$. \qed
\begin{theorem} \label{t-main3}
For every $f\in L^2(Q)$ there exists a unique solution $v_0\in C\bigl(\Gamma,H^1_{per,0}(Q)
\bigr)$ of (\ref{eq-operator-0}) and (\ref{var-5}). The mapping $f\mapsto v_0$ is bounded.
\end{theorem}
\textbf{Proof:} As in the proof of Lemma~\ref{l-main2} we introduce the operator 
$T_0:L^2(Q)\to C\bigl(\Gamma,H^1_{per,0}(Q)\bigr)$ by $T_0g=w_g$ where $w_g\in 
C\bigl(\Gamma,H^1_{per,0}(Q)\bigr)$ is the unique solution of 
$$ (I-K_{0,\alpha})w_g(\cdot,\alpha)\ =\ R_\alpha g\quad\mbox{in }H^1_{per,0}(Q)\mbox{ for 
all }\alpha\in\Gamma\,. $$
The existence and boundedness of $T_0$ is again assured by Lemma~\ref{l-main1}. Then we 
write (\ref{eq-operator-0}) as the fixpoint equation
$$ v_0\ =\ T_0(f+k^2qPv_0)\ =\ T_0f\ +\ k^2T_0(qPv_0)\quad\mbox{in }
C\bigl(\Gamma,H^1_{per,0}(Q)\bigr) $$
which is again of Fredholm type because of the compactness of $P$. Therefore, it suffices to 
prove uniqueness. Let $v_0\in C\bigl(\Gamma,H^1_{per,0}(Q)\bigr)$ be a solution corresponding 
to $f=0$. Define $u$ by (\ref{u-def}) in $W_{h_0}$. Then $u$ satisfies $\Delta u+k^2nu=-k^2qu$
and the open waveguide radiation condition. A well known uniqueness result of 
Furuya (\cite{Furuya2020}) yields that $u$ vanishes identically. This implies that $Pv_0=u$ 
vanishes and thus $v_0=0$ by the injectivity of $I-K_{0,\alpha}$. \qed

\section{The PML-Method}

We define the PML-operator as in \cite{Zhang2021b} by choosing a complex-valued function 
$\hat{s}\in C^1[0,h_0+\tau]$ for some $\tau>0$ (thickness of the PML-layer) such that 
$\hat{s}(x_2)=0$ for $x_2\leq h_0$ and $\Re\hat{s}(x_2)>0$ and $\Im\hat{s}(x_2)>0$ for 
$x_2\in[h_0,h_0+\tau]$. 
 With such a function $\hat{s}$ and parameter $\rho>0$ we define the 
function $s_\rho$ by $s_\rho(x_2)=1+\rho\hat{s}(x_2)$ for $x_2\in[0,h_0+\tau]$. In the 
following we choose the particular function
$$ s_\rho(x_2)\ =\ 1+\,\rho\,e^{i\pi/4}\biggl(\frac{x_2-h_0}{\tau}\biggr)^m\mbox{ for }
x_2\in[h_0,h_0+\tau],\ s_\rho(x_2)\ =\ 1\mbox{ for }x_2\leq h_0\,, $$
where $m\geq 1$ is some integer and $\rho>0$ is the PML-parameter which is assumed to be 
large. With this function $s_\rho$ we define the operator $\Delta_\rho$ by
$$ \Delta_\rho u\ :=\ \frac{\partial^2u}{\partial x_1^2}\ +\ \frac{1}{s_\rho}\frac{\partial u}
{\partial x_2}\biggl(\frac{1}{s_\rho}\frac{\partial}{\partial x_2}\biggr)\,, $$
and look at the boundary value problem to determine $u=u_{\rho}\in H^1_{loc}(W_{h_0+\tau})$ 
with
\begin{equation} \label{eqn-PML}
\Delta_\rho u+k^2(n+q)u\ =\ -f\mbox{ in }W_{h_0+\tau}\,,\quad u=0\mbox{ for }
x\in\partial W_{h_0+\tau}\,.
\end{equation}
Therefore, instead of \eqref{eq_de:a}--\eqref{eq_de:c} we consider
\begin{eqngroup}\begin{eqnarray}
& & \Delta_\rho v(\cdot,\alpha,\rho)\ +\ 2i\alpha\,\partial_1v(\cdot,\alpha,\rho)\ +\
[k^2\,n-\alpha^2]\,v(\cdot,\alpha,\rho) \nonumber \\ 
& = & -f\,e^{-i\alpha x_1}\ -\ k^2\,q\,e^{-i\alpha x_1}\int\limits_\Gamma
v(\cdot,\beta,\rho)\,e^{i\beta x_1}\,d\beta\quad\mbox{in }(0,2\pi)\times(0,h_0+\tau)\,, 
\label{eq_PML:a} \\
& & v(\cdot,\alpha,\rho)\ =\ 0\quad\mbox{for $x_2=0$ or $x_2=h_0+\tau$}\,, \label{eq_PML:b}
\end{eqnarray} \end{eqngroup}
for all $\alpha\in\Gamma$. To compare it with the exact solution we transform this problem to 
a problem on $Q$ with the help of the Dirichlet to Neumann map. It is well known that the 
Dirichlet to Neumann map $\Lambda_\alpha:H^{1/2}_{per}(\gamma)\to 
H^{-1/2}_{per}(\gamma)$ for the original problem (where $\gamma=(0,2\pi)\times\{h_0\}$), 
defined as
\begin{equation} \label{DtN}
(\Lambda_\alpha\phi)(x_1,h_0)\ :=\ \frac{i}{\sqrt{2\pi}}\sum_{\ell\in\ganz}
\sqrt{k^2-(\ell+\alpha)^2}\,\phi_\ell(h_0)\,e^{i\ell x_1}\,,
\end{equation}
is replaced by
\begin{equation} \label{DtN-PML}
(\Lambda_{\alpha,\rho}\phi)(x_1,h_0)\ :=\ \frac{i}{\sqrt{2\pi}}\sum_{\ell\in\ganz}
\sqrt{k^2-(\ell+\alpha)^2}\,\coth\bigl(-i\sqrt{k^2-(\ell+\alpha)^2}\sigma_\rho\bigr)\,
\phi_\ell(h_0)\,e^{i\ell x_1}\,,
\end{equation}
where $\sigma_\rho=\int_{h_0}^{h_0+\tau}s_\rho(t)dt=\tau+\rho\chi$ with 
$\chi=(1+i)\frac{\tau}{m+1}$. Therefore, \eqref{eq_PML:a}, \eqref{eq_PML:b} is equivalent 
to (compare with \eqref{eq_de:a}--\eqref{eq_de:c})
\begin{eqngroup}\begin{eqnarray}
& & \Delta v(\cdot,\alpha,\rho)\ +\ 2i\alpha\,\partial_1v(\cdot,\alpha,\rho)\ +\
[k^2\,n-\alpha^2]\,v(\cdot,\alpha,\rho) \nonumber \\ 
& = & -f\,e^{-i\alpha x_1}\ -\ k^2\,q\,e^{-i\alpha x_1}\int\limits_\Gamma
v(\cdot,\beta,\rho)\,e^{i\beta x_1}\,d\beta\quad\mbox{in }Q\,, \label{eq_de_PML:a} \\
v(\cdot,\alpha,\rho) & = & 0\quad\mbox{for }x_2=0\,, \label{eq_de_PML:b} \\
\partial_2v_\ell(h_0,\alpha,\rho) & = & i\sqrt{k^2-(\ell+\alpha)^2}\,
\coth\bigl(-i\sqrt{k^2-(\ell+\alpha)^2}\sigma_\rho\bigr)\,v_\ell(h_0,\alpha,\rho)\,, 
\label{eq_de_PML:c}
\end{eqnarray} \end{eqngroup}
for all $\alpha\in\Gamma$. Its variational form is given by
\begin{eqnarray}
& & \int\limits_Q\biggl[\nabla v(\cdot,\alpha,\rho)\cdot\nabla\overline{\psi}-2i\alpha\,
\partial_1v(\cdot,\alpha,\rho)\,\overline{\psi}-[k^2\,n-\alpha^2]\,
v(\cdot,\alpha,\rho)\,\overline{\psi}\biggr]\,dx\ - \label{var-6} \\
& &  -\ i\sum_{\ell\in\ganz}\sqrt{k^2-(\ell+\alpha)^2}\,
\coth\bigl(-i\sqrt{k^2-(\ell+\alpha)^2}\sigma_\rho\bigr)v_\ell(h_0,\alpha,\rho)\,
\overline{\psi_\ell(h_0)} \nonumber \\ 
& = & \int\limits_Qf\,e^{-i\alpha x_1}\,\overline{\psi}\,dx\ +\ k^2
\int\limits_Q \,q\,e^{-i\alpha x_1}\,\overline{\psi}\int\limits_\Gamma
v(\cdot,\beta,\rho)\,e^{i\beta x_1}\,d\beta\,dx \nonumber
\end{eqnarray}
for all $\psi\in H^1_{per,0}(Q)$ and $\alpha\in\Gamma$. First we show
\begin{lemma} \label{l-diff}
There exist $c>0$ and $\mu>0$ with
\begin{equation} \label{eq-diff}
\Vert\Lambda_{\alpha,\rho}-\Lambda_\alpha\Vert_{H^{1/2}(\gamma)\to H^{-1/2}(\gamma)}\ \leq\
c\,e^{-\mu\rho}\quad\mbox{for all $\rho\geq 0$ and }\alpha\in\Gamma\,.
\end{equation}
\end{lemma}
\textbf{Proof:} We observe that the difference $\Lambda_{\alpha,\rho}-\Lambda_\alpha$ conatains
\begin{eqnarray*}
& & \sqrt{k^2-(\ell+\alpha)^2}\,\bigl[\coth\bigl(-i\sqrt{k^2-(\ell+\alpha)^2}\sigma_\rho\bigr)
-1\bigr] \\
& = & \frac{2\sqrt{k^2-(\ell+\alpha)^2}}{e^{-2i\sqrt{k^2-(\ell+\alpha)^2}\sigma_\rho}-1}\ =\
\frac{2\sqrt{k^2-(\ell+\alpha)^2}}{e^{-2i\sqrt{k^2-(\ell+\alpha)^2}(\tau+\rho\chi)}-1}\,.
\end{eqnarray*}
We set $\omega:=\ell+\alpha=\omega_1+i\omega_2$ with $\omega_j\in\real$ for abbreviation and 
$t(\omega)=2\sqrt{k^2-\omega^2}=t_1(\omega)+it_2(\omega)$ (with $t_j(\omega)\in\real$) and 
note that, by the choice of $\Gamma$, there exists $c_1>0$ with $\bigl||\omega|-k\bigr|
\geq c_1$ for all $\omega:=\ell+\alpha$, $\ell\in\ganz$, $\alpha\in\Gamma$. Then, with 
$\chi_1=\Re\chi=\Im\chi=\frac{\tau}{m+1}$,
$$ \frac{2\sqrt{k^2-(\ell+\alpha)^2}}{e^{-2i\sqrt{k^2-(\ell+\alpha)^2}(\tau+\rho\chi)}-1}\ =\
\frac{t(\omega)}{e^{-it(\omega)(\tau+\rho\chi)}}\ =\ \frac{t(\omega)}
{e^{-i(t_1\tau+t_1\rho\chi_1-t_2\rho\chi_1)}e^{\rho\chi_1(t_1+t_2)+t_2\tau}-1} $$
and thus
$$ \left|\frac{2\sqrt{k^2-(\ell+\alpha)^2}}{e^{-2i\sqrt{k^2-(\ell+\alpha)^2}(\tau+\rho\chi)}-1}
\right|\ \leq\ \frac{|t(\omega)|}{e^{\rho\chi_1(t_1+t_2)+t_2\tau}-1} $$
provided $t_1+t_2>0$. We consider two cases:
\smalf
(a) $|\omega|\geq k+c_1$. Then $\omega_1^2+\omega_2^2\geq k^2+c_1^2$ and thus $t_1^2-t_2^2=
\Re(t^2)=4\Re(k^2-\omega^2)=4[k^2-\omega_1^2+\omega_2^2]\leq 4\omega_2^2-4c_1^2\leq-3c_1^2
=:-c_2$ provided $\delta$ (in the definition of $\Gamma$) is small enough. Therefore, 
$t_2>0$ (by the choice of the square root function) and $|t_1|\leq t_2-c_3$ where $c_3>0$ 
is independent of $\omega$. Therefore,
$$ \frac{|t(\omega)|}{e^{\rho\chi_1(t_1+t_2)+t_2\tau}-1}\ \leq\
\frac{|t(\omega)|}{e^{\rho\chi_1(t_2-(t_2-c_3))}-1}\ =\ 
\frac{|t(\omega)|}{e^{\rho\chi_1c_3}-1}\ \leq\ ce^{-\rho\chi_1c_3/2}\,. $$
(b) $|\omega|\leq k+c_1$. Then $\omega_1^2+\omega_2^2\leq k^2-c_1^2$ and thus $t_1^2-t_2^2=
\Re(t^2)=4\Re(k^2-\omega^2)=4[k^2-\omega_1^2+\omega_2^2]\geq 4\omega_2^2+4c_1^2\geq 4c_1^2
=:c_4$. Therefore, $t_1>0$ (again by the choice of the square root function) and 
$|t_2|\leq t_1-c_5$ where $c_5>0$ is independent of $\omega$. Therefore,
$$ \frac{|t(\omega)|}{e^{\rho\chi_1(t_1+t_2)+t_2\tau}-1}\ \leq\
\frac{|t(\omega)|}{e^{\rho\chi_1(t_1-(t_1-c_5))+t_2\tau}-1}\ =\ 
\frac{|t(\omega)|}{e^{\rho\chi_1c_5+t_2\tau}-1}\ \leq\ ce^{-\rho\chi_1c_5/2} $$
because $t_2=t_2(\omega)$ is bounded in this case. \qed
\begin{theorem}
There exists $\rho_0>0$ such that (\ref{var-6}) has a unique solution $v(\cdot,\cdot,\rho)
\in C\bigl(\Gamma,H^1_{per}(Q)\bigr)$ for all $\rho\geq\rho_0$. Furthermore, there exist 
$c,\mu>0$ with 
\begin{equation}
\max\limits_{\alpha\in\Gamma}\Vert v(\cdot,\alpha,\rho)-v_0(\cdot,\alpha)\Vert_{H^1(Q)}\
\leq\ c\,e^{-\mu\rho}\,,\quad\rho\geq\rho_0\,.
\end{equation}
\end{theorem}
As a corollary we transform this result into the fields $u$. Let
\begin{eqnarray*}
u(x_1,x_2) & = & \int\limits_\Gamma v_0(x,\alpha)\,e^{i\alpha x_1}\,\,d\alpha\,, \\
u(x_1,x_2,\rho) & = & \int\limits_\Gamma v(x,\alpha,\rho)\,e^{i\alpha x_1}\,\,d\alpha
\end{eqnarray*}
for $x\in W_{h_0}$. Then we have
\begin{corollary} 
There exists $\rho_0>0$ such that for all $R>0$ there exist $c=c(R)$ and $\mu=\mu(R)>0$ 
such that
\begin{equation}
\Vert u(\cdot,\rho)-u\Vert_{H^1(Q_R)}\ \leq c\,e^{-\mu\rho}\,,\quad\rho\geq\rho_0\,.
\end{equation}
Here, $Q_R=(-R,R)\times(0,h_0)$.
\end{corollary}

\section{Numerical Implementation}

\subsection{A Hybrid Spectral-FD-Method}

We use a hybrid spectral-FD method for the approximation of \eqref{eq_PML:a}, 
\eqref{eq_PML:b}; that is, we expand $v(\cdot,\alpha)$ into the form
$$ v(x_1,x_2;\alpha)\ =\ \sum_{\ell=-\infty}^\infty v_\ell(x_2,\alpha)\,e^{i\ell x_1}\,. $$
We expand $n(x_1,x_2)=\sum_{m=-\infty}^\infty n_m(x_2)\,e^{im x_1}$ and 
$f(x_1,x_2)e^{-i\alpha x_1}=\sum_{\ell=-\infty}^\infty f_\ell(x_2,\alpha)\,e^{i\ell x_1}$ 
and $q(x_1,x_2)e^{-i(\alpha-\beta)x_1}=\sum_{m=-\infty}^\infty 
q_m(x_2,\alpha-\beta)\,e^{im x_1}$ we observe that
$v_\ell(x_2,\alpha)$ satisfies the following coupled system of ordinary differential 
equations
\begin{eqnarray*}
& & a(x_2)\,v_\ell^{\prime\prime}(x_2,\alpha)\ +\ b(x_2)\,v_\ell^\prime(x_2,\alpha)
-(\ell+\alpha)^2v_\ell(x_2,\alpha)+
k^2\sum_{m=-\infty}^\infty v_m(x_2,\alpha)\,n_{\ell-m}(x_2) \\ 
& = & f_\ell(x_2,\alpha)\ -\ k^2\sum_{m=-\infty}^\infty\int\limits_\Gamma 
v_m(x_2,\beta)\,q_{\ell-m}(x_2,\alpha-\beta)\,d\beta\,,\quad x_2\in(0,h_0+\tau)\,,
\end{eqnarray*}
and $v_\ell(0,\alpha)=v_\ell(h_0+\tau,\alpha)=0$ for all $\ell\in\ganz$ and 
$\alpha\in\Gamma$. Here we have set $a=\frac{1}{s^2}$ and 
$b=\frac{1}{s}\left(\frac{1}{s}\right)^\prime$.
\smalf
We restrict $\ell$ to $\ell\in\{-L,\ldots,L\}$ for some 
$L\in\nat$, truncate the series to $\sum_{m=-L}^L$, and approximate the integral 
by a quadrature rule $\int_\Gamma g(\beta)\,d\beta\approx\sum_{\mu=0}^M 
w_\mu\,g(\gamma_\mu)$ for some $\gamma_\mu\in\Gamma$ and weights $w_\mu$. Setting 
$\alpha=\gamma_\nu$ and $v_\ell(x_2,\gamma_\nu)=v_{\ell,\nu}(x_2)$ the system now 
reduces to
\begin{eqnarray*}
& & a(x_2)\,v_{\ell,\nu}^{\prime\prime}(x_2)\ +\ b(x_2)\,v_{\ell,\nu}^\prime(x_2)
-(\ell+\gamma_\nu)^2v_{\ell,\nu}(x_2)+
k^2\sum_{m=-L}^L v_{m,\nu}(x_2)\,n_{\ell-m}(x_2) \\ 
& = & f_\ell(x_2,\gamma_\nu)\ -\ k^2\sum_{m=-L}^L\sum_{\mu=0}^M w_\mu\, 
v_{m,\mu}(x_2)\,q_{\ell-m}(x_2,\gamma_\nu-\gamma_\mu)\,,\quad x_2\in(0,h_0+\tau)\,,
\end{eqnarray*}
and $v_{\ell,\nu}(0)=v_{\ell,\nu}(h_0+\tau)=0$ for all $\ell=-L,\ldots,L$ and 
$\nu=0,\ldots,M$. 
\medlf
This coupled system of boundary value problems for ordinary differential equations is 
approximately solved by a finite difference method: Let $n_y\in\nat$ and 
$h_y=(h_0+\tau)/n_y$. Then $v_{\ell,\nu,j}\approx v_\ell(jh_y,\gamma_\nu)$ is 
determined by
\begin{eqnarray}
& & \frac{a(jh_y)}{h_y^2}\bigl[v_{\ell,\nu,j+1}+v_{\ell,\nu,j-1}-2v_{\ell,\nu,j}\bigr] 
+ \frac{b(jh_y)}{2h_y}\bigl[v_{\ell,\nu,j+1}-v_{\ell,\nu,j-1}\bigr] \label{eq:system} \\ 
& & +\ k^2\sum_{m=-L}^L v_{m,\nu,j}\,n_{\ell-m}(jh_y)\ -\
(\gamma_\nu+\ell)^2\,v_{\ell,\nu,j} \nonumber \\
& = & -f_\ell\bigl(jh_y,\gamma_\nu)\ -\ 
k^2\sum_{\mu=0}^Mw_\mu\sum_{m=-L}^L q_{\ell-m}(jh_y,\gamma_\nu-\gamma_\mu)\,
v_{m,\mu,j} \nonumber
\end{eqnarray}
for $\ell=-L,\ldots,L$ and $j=1,\ldots,n_y-1$ and $\nu=0,\ldots,M$. 
\smalf
If $n$ and the perturbation $q$ are separable; that is, of the forms
$n(x)=n^{(1)}(x_1)n^{(2)}(x_2)$ and $q(x)=q^{(1)}(x_1)q^{(2)}(x_2)$, respectively, then
the sum on the left hand side is, for fixed $\nu$ and $j$, just a matrix-vector 
multiplication of $P:=(n^{(1)}_{\ell-m})_{\ell,m=1}^L$ with $(v_{m,\nu,j})_{m=1}^L$. 
Furthermore, we set 
$$ Q_{(\ell,\nu),(m,\mu)}\ :=\ q^{(1)}_{\ell-m}(\gamma_\nu-\gamma_\mu)\,w_\mu\quad
\mbox{and}\quad v_{(m,\mu)}:=v_{m,j,\mu} $$
for fixed $j$. Re-arranging the pairs $(\ell,\nu)$ and $(m,\mu)$ into one index each 
the double sum on the right hand side is just a matrix-vector multiplication. The matrices
$P$ and $Q$ can be computed beforehand.
\biglf
\begin{remark} \label{rem:D2N}
Usually, the PML method is used to replace the non-local boundary
condition $\partial_2 u=\Lambda_\alpha u$ with the Dirichlet-to-Neumann map 
$\Lambda_\alpha$ on $(0,2\pi)\times \{h_0+\tau\}$ by local conditions in the PML-layer
$(0,2\pi)\times (h_0,h_0+\tau)$. For the spectral method -- which is already non-local 
with respect to $x_1$ -- the condition $\partial_2 u=\Lambda_\alpha u$ is local and given
by
$$ \partial_{x_2}v_m(h_0+\tau,\alpha)=i\sqrt{k^2-(m+\alpha)^2}\,
v_m(h_0+\tau,\alpha)\,. $$
Therefore, we set $a\equiv 1$ and $b\equiv 0$ in \eqref{eq:system} and discretize the 
Rayleigh condition as
$$ \frac{1}{2h_y}\bigl[v_{\ell,\nu,n_y}-v_{\ell,\nu,n_y-2}\bigr]\ =\ 
\sqrt{k^2-(\ell+\gamma_\nu)^2}\,v_{\ell,\nu,n_y-1}\,. $$
We express $v_{\ell,\nu,n_y}$ by $v_{\ell,\nu,n_y-1}$ and $v_{\ell,\nu,n_y-2}$ and 
substitute this into equation \eqref{eq:system} for $j=n_y-1$. This yields again a system 
of $n_y-1$ unknowns $v_{\ell,\nu,j}$, $j=1,\ldots,n_y-1$.
\end{remark}

\subsection{A Particular Example}

\begin{figure}[h] 
\includegraphics[width=0.5\textwidth]{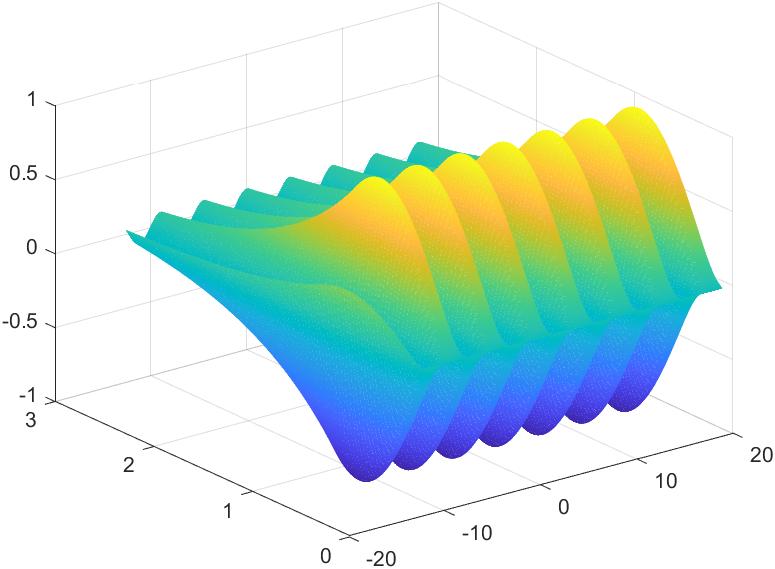}\hfill
\includegraphics[width=0.4\textwidth]{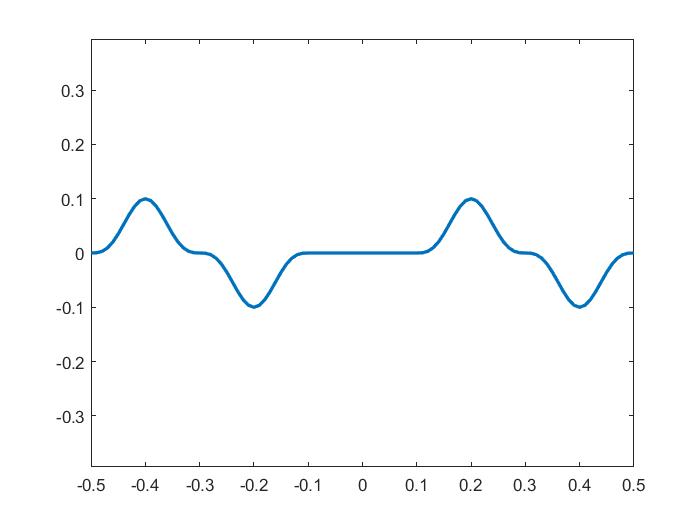}
\caption{\label{fig:mode} The real part of the mode $\phi^+$ (left) and the 
curve $\Gamma$ (right)}
\end{figure}

\begin{figure} 
\includegraphics[width=0.4\textwidth]{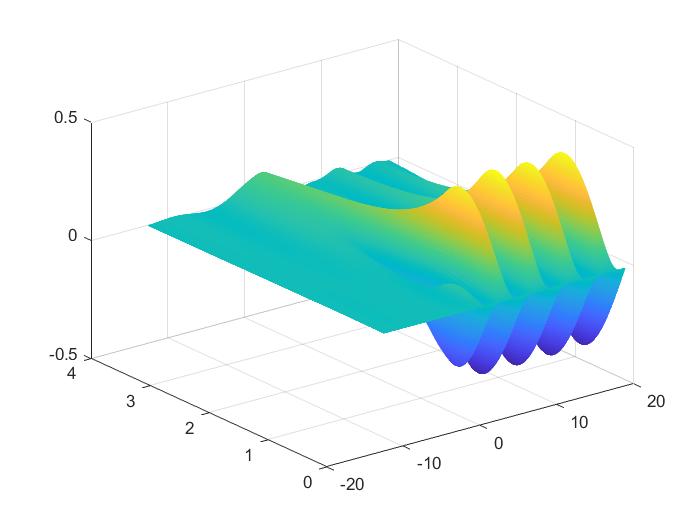}
\includegraphics[width=0.4\textwidth]{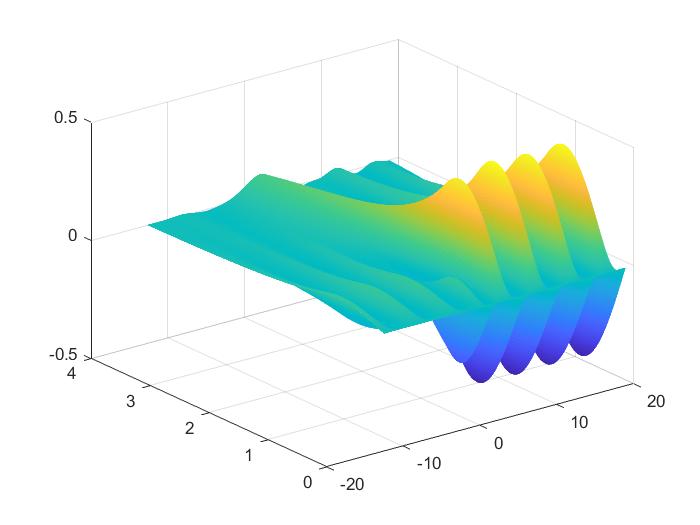}

\includegraphics[width=0.4\textwidth]{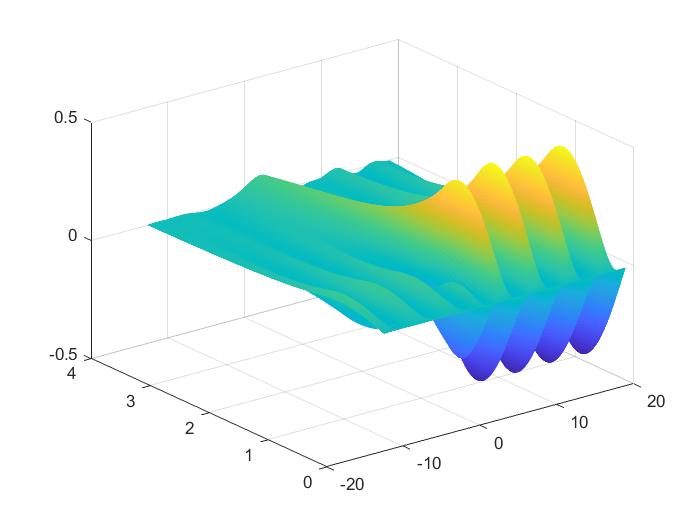}
\includegraphics[width=0.4\textwidth]{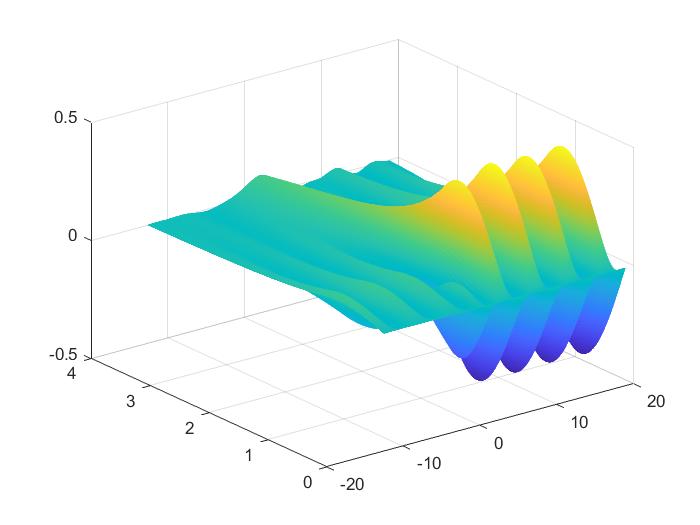}
\caption{\label{fig:iter} Real parts of $u^{(1)}$ (i.e. Born), $u^{(2)}$, $u^{(5)}$, 
and $u^{(9)}$.}
\end{figure}

\begin{figure}[h] 
\includegraphics[width=0.4\textwidth]{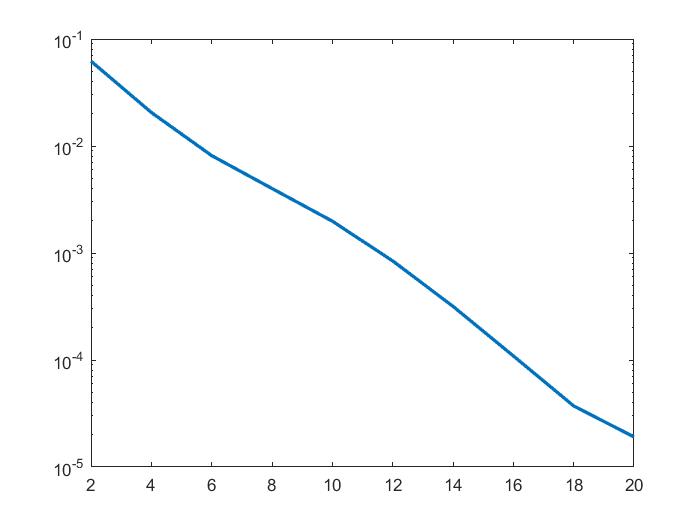}
\caption{\label{fig:error} Error as a function of $\rho$ between the 4th iterates of
the PML-method and the D2N-boundary condition}
\end{figure}

We consider a simple example of the scattering of a mode corresponding to a constant 
refractive index $n>1$ in the layer $\real\times(0,1)$ by a local perturbation $q$.
The propagating modes $\phi^\pm$ are assumed to be quasi-periodic solutions of 
$\Delta\phi+k^2n\phi=0$ in $\real\times(0,1)$, $\Delta\phi+k^2\phi=0$ for $x_2>1$,
satisfy homogeneous boundary conditions for $x_2=0$ and the Rayleigh expansion for 
$x_2>1$ and continuity conditions for $\phi$ and $\partial_2\phi$ for $x_2=1$. They
are of the form
$$ \phi^\pm(x)\ =\ e^{\pm i\omega x_1}\left\{\begin{array}{cl} \sin\sqrt{nk^2-\omega^2}\,
e^{-\sqrt{\omega^2-k^2}\,(x_2-1)}\,, & x_2>1\,, \\
\sin\bigl(\sqrt{nk^2-\omega^2}\,x_2\bigr)\,, & 0<x_2<1\,. \end{array}\right. $$
Here, $\omega\in(k,\sqrt{n}k)$ is a propagative wave number if $k$, $n$, and 
$\omega$ satisfy the equation
$$ \sqrt{\omega^2-k^2}\,\sin\bigl(\sqrt{nk^2-\omega^2}\bigr)\ +\ \sqrt{nk^2-\omega^2}\,
\cos\bigl(\sqrt{nk^2-\omega^2}\bigr)\ =\ 0\,. $$
We choose the wavenumber $k>0$ and frequency $\omega>k$ arbitrarily and determine 
the constant $n$ such that $\sqrt{k^2n-\omega^2}\cot\sqrt{k^2n-\omega^2}=
-\sqrt{\omega^2-k^2}$. Therefore, we first determine $z\in(\pi/2,\pi)$ such that 
$z\cot z=-\sqrt{\omega^2-k^2}$ and then determine $n$ from $z=\sqrt{k^2n-\omega^2}$;
that is, $n=\frac{z^2+\omega^2}{k^2}$. We consider the 
particular example $k=0.8$ and $\omega=1.4$. Then $n\approx 9.8$. We observe that
$\hat\alpha=0.4$ and $\kappa=-0.2$ for this example. $\phi^+$ is right-going, $\phi^-$ is 
left-going. The real part of $\phi^+$ is shown in Figure~\ref{fig:mode} (left).
\smalf
We choose $h_0=2.5$ and the curve $\Gamma$ a bit different than in the paper. The half circle 
$\{z\in\cmplx:|z-a|=\delta,\ \pm\Im z<0\}$ for $a\in\real$ is replaced by the curve 
$\{t\mp i\tilde\gamma(t)\in\cmplx:t\in(a-\delta,a+\delta)\}$ where 
$$ \tilde\gamma(t)\ =\ \eps\sin^n\biggl[\frac{\pi}{2}\biggl(\frac{t-a}{\delta}+1
\biggr)\biggr]\,,\quad a-\delta<t<a+\delta\,. $$
For our example we take $n=3$ and $\eps=\delta=0.1$. This leads to a global 
parameterization $\alpha=\gamma(t)$, $t\in(0,T)$, of $\Gamma$. For 
the numerical integration of $T-$periodic functions on $\Gamma$ we use the trapezoidal 
rule; that is, $\int_\Gamma g(\alpha)\,d\alpha\approx \sum_{\mu=0}^Mw_\mu\,g(\gamma_\mu)$ 
with $\gamma_\mu=\gamma(\mu T/M)$ for $\mu=0,\ldots,M$ and $w_\mu=
\frac{T}{M}\gamma^\prime(\mu T/M)$ for $\mu=1,\ldots,M-1$ and 
$w_0=w_M=\frac{T}{2M}\gamma^\prime(0)=\frac{T}{2M}$.
\smalf
The curve $\Gamma$ is shown in Figure~\ref{fig:mode} (right).
\medlf
\medlf
We consider the scattering of the mode $\phi^+$ by a local perturbation $q=q(x)$. 
The total field $u^{tot}=u+\phi^+$ as the sum of the incident field $u^{inc}=\phi^+$,
and the scattered field $u$ satisfies $\Delta u^{tot}+k^2(n+q)u^{tot}=0$. Therefore, 
the scattered field $u$ satisfies $\Delta u+k^2(n+q)u=-k^2q\,\phi^+$ and the 
radiation condition. In this case the source is given by $f=k^2q\,\phi^+$. Writing 
the previous equation as 
$$ \Delta u+k^2nu\ =\ -k^2q\,\phi^+\ -\ k^2q\,u $$
we can -- under a smallness assumption on $k^2q$ -- iterate this fixpoint 
equation for $u$. The first iteration is the Born approximation $u^B$, defined 
as the solution of $\Delta u^B+k^2nu^B=-k^2q\,\phi^+$. As a particular example
for the perturbation $q$ we choose $q(x_1,x_2)=q_0\sin^2(2\omega x_1)$ for $(x_1,x_2)\in 
(0,\pi/\omega)\times(0.2,0.7)$ and zero else. Here, $q_0>0$ is a parameter. This 
perturbation has the property that
\begin{eqnarray}
\int\limits_Qf(x)\,\overline{\phi^-(x)}\,dx & = & k^2\int\limits_Qq(x)\,
\phi^+(x)\,\overline{\phi^-(x)}\,dx \nonumber \\ 
& = & k^2q_0\int\limits_0^{\pi/\omega}\sin^2(2\omega x_1)\,e^{2i\omega x_1}dx_1\,
\int\limits_{0.2}^{0.7}\sin^2(zx_2)\dx_2\ =\ 0\quad\mbox{and} \label{eq:ortho1} \\
\int\limits_Qf(x)\,\overline{\phi^+(x)}\,dx & = & k^2\int\limits_Qq(x)\,
|\phi^+(x)|^2\,dx\ >\ 0\,. \nonumber
\end{eqnarray}
From \eqref{eq:coeff_a} we observe that the Born approximation $u^B$ has a 
right-going mode but no left-going mode.
\biglf
We choose the parameter $q_0=n/2$; that is, the perturbation is (in the $\infty-$norm)
of order 50\%. For constant $n$ we rewrite the iterative scheme of \eqref{eq:system} as
\begin{eqnarray}
& & \frac{a(jh_y)}{h_y^2}\bigl[v^{(t+1)}_{\ell,\nu,j+1}+v^{(t+1)}_{\ell,\nu,j-1}-
2v^{(t+1)}_{\ell,\nu,j}\bigr] 
+ \frac{b(jh_y)}{2h_y}\bigl[v^{(t+1)}_{\ell,\nu,j+1}-v^{(t+1)}_{\ell,\nu,j-1}\bigr] 
\nonumber \\ 
& & +\ [k^2n-(\gamma_\nu+\ell)^2]\,v^{(t+1)}_{\ell,\nu,j} \label{eq:system2} \\
& = & -f_\ell\bigl(jh_y,\gamma_\nu)\ -\ 
k^2\sum_{\mu=0}^Mw_\mu\sum_{m=-L}^L q_{\ell-m}(jh_y,\gamma_\nu-\gamma_\mu)\,
v^{(t)}_{m,\mu,j} \nonumber
\end{eqnarray}
for $\ell=-L,\ldots,L$ and $j=1,\ldots,n_y-1$ and $\nu=0,\ldots,M$, and $t=0,1,\ldots.$.
Setting $v^{(0)}_{m,\mu,j}=0$ we obtain $v^{(1)}_{m,\mu,j}$ as the Born approximation.
The field $u^{(t)}$ is then computed as the integral 
$$ u^{(t)}(x_1,x_2)\ =\ \sum_{\ell=-L}^L\int_\Gamma v^{(t)}_\ell(x_2,\alpha)\,
e^{i(\ell+\alpha)x_1}d\alpha\ \approx\ \sum_{\ell=-L}^L\sum_{\mu=0}^M 
v^{(t)}_\ell(x_2,\gamma_\mu)\,w_\mu\,e^{i(\ell+\gamma_\mu)x_1}\,. $$
\medlf
As a comparison with the PML-method for different parameters $\rho$ we first 
use the Dirichlet-to-Neumann map as a boundary condition at $x_2=h_0+\tau$ as explained 
in Remark~\ref{rem:D2N}. We computed the relative error 
$e_t:=\Vert u^{(t+1)}-u^{(t)}\Vert_\infty/\Vert u^{(t+1)}\Vert_\infty$ on
$(-4\pi,6\pi)\times(0,4)$ between consecutive iterates for $t=1,\ldots,8$ as
$$ \begin{array}{ccccccccc}
t: & 1 & 2 & 3 & 4 & 5 & 6 & 7 & 8 \\
e_t: & 0.2276 & 0.0674 & 0.0205 & 0.0048 & 0.0017 & 0.0005 & 0.0001 & 0.0000 
\end{array} $$
and show the iterations (real parts) $u^{(1)}$ (i.e. Born), $u^{(2)}$, $u^{(5)}$, and 
$u^{(9)}$ on
$(-4\pi,6\pi)\times(0,4)$ for the Dirichlet-to-Neumann map in Figure~\ref{fig:iter}.
\smalf
We clearly observe that the Born approximation; that is, the first iteration, has no 
left going mode because $f$ is orthogonal to $\phi^-$, see \eqref{eq:ortho1}. The further 
iterations produce right hand sides which are not orthogonal to $\phi^-$ anymore and,
therefore, left going modes of small amplitudes appear.
\medlf
Next, we implemented the PML-method of \eqref{eq:system2} to show the dependence of the 
result with respect to the parameter $\rho$.. We iterated this system and compared the 
4th iteration of the PML-method for parameter $\rho=2,4,6\ldots,20$ with the 4th 
iteration of the system with Dirichlet-to-Neumann boundary condition. The 
(semi-logarithmic) plot of the relative errors 
$$ \max\{|u^{(4)PML}(x)-u^{(4)D2N}(x)|:x\in[-4\pi,6\pi]\times[0,0.9*2.5]\}\ /\
\Vert u^{(4)D2N}\Vert_\infty $$
of the 4th iterations is shown in Figure~\ref{fig:error}. One clearly observes the exponential
decay.
Finally, we summarize the settings of the parameters in our simulations:
$$ \begin{array}{lc}
k \mbox{ (wavenumber):} & 0.8 \\
\omega \mbox{ (propagative wavenumber):} & 1.4 \\
n \mbox{ (refractive index, determined from $k$ and $\omega$):} & \approx 9.8 \\
h_0 \mbox{ (height of domain without PML-layer):} & 2.5 \\
\tau \mbox{ (height of PML-layer):} & 1.5 \\
s_\rho \mbox{ (PML-function):} & 1+\,\rho\,e^{i\pi/4}(x_2-h_0)^3/\tau^3 \\
n_y \mbox{ (number of discretization points in $(0,h_0+\tau)$):} & 512 \\
L \mbox{ (summation bound for spectral representation):} & 7 \\
M \mbox{ (number of discretization/integration points on $\Gamma$):} & 101 \\
\end{array} $$
We observed that the results do not depend crucially on the values of $n_y$, 
$L$, or $M$, or the parameters of $s_\rho$.

\bibliographystyle{plain}
\bibliography{/Users/zhang/Documents/desktop/my_papers/KZ_PML_layer/ip-biblio.bib}

\providecommand{\noopsort}[1]{}
\begin{thebibliography}{10}

\bibitem{Beren1994}
J.-P. Berenger.
\newblock A perfectly matched layer for the absorption of electromagnetic
  waves.
\newblock {\em Journal of Computational Physics}, 114(2):185–200, 1994.

\bibitem{Chand2009}
S.~N. Chandler-Wilde and P.~Monk.
\newblock The {P}{M}{L} for rough surface scattering.
\newblock {\em Applied Numerical Mathematics}, 59:2131--2154, 2009.

\bibitem{Chen2003}
Z.~Chen and H.~Wu.
\newblock An adaptive finite element method with perfectly matched absorbing
  layers for the wave scattering by periodic structures.
\newblock {\em SIAM Journal on Numerical Analysis}, 41(3):799–826, 2003.

\bibitem{Coatl2012}
J.~Coatl{\'e}ven.
\newblock {Helmholtz equation in periodic media with a line defect}.
\newblock {\em {J. Comp. Phys.}}, 231:1675--1704, 2012.

\bibitem{Colto2019}
D.~L. Colton and R.~Kress.
\newblock {\em Inverse acoustic and electromagnetic scattering theory}.
\newblock Springer, 4. edition, 2093.

\bibitem{Ehrha2009b}
M.~Ehrhardt, H.~Han, and C.~Zheng.
\newblock Numerical simulation of waves in periodic structures.
\newblock {\em Commun. Comput. Phys.}, 5:849--870, 2009.

\bibitem{Ehrha2009a}
M.~Ehrhardt, J.~Sun, and C.~Zheng.
\newblock Evaluation of scattering operators for semi-infinite periodic arrays.
\newblock {\em Commun. Math. Sci.}, 7(2):347--364, 2009.

\bibitem{Ehrha2008a}
M.~Ehrhardt and C.~Zheng.
\newblock Exact artificial boundary conditions for problems with periodic
  structures.
\newblock {\em J. Comput. Phys.}, 227(14):6877--6894, 2008.

\bibitem{Fliss2013}
S.~Fliss.
\newblock A dirichlet-to-neumann approach for the exact computation of guided
  modes in photonic crystal waveguides.
\newblock {\em SIAM Journal on Scientific Computing}, 35(2):B438–B461, 2013.

\bibitem{Fliss2009a}
S.~Fliss and P.~Joly.
\newblock Exact boundary conditions for time-harmonic wave propagation in
  locally perturbed periodic media.
\newblock {\em Appl. Numer. Math.}, 59:2155--2178, 2009.

\bibitem{Fliss2015}
S.~Fliss and P.~Joly.
\newblock {Solutions of the time-harmonic wave equation in periodic waveguides:
  asymptotic behaviour and radiation condition}.
\newblock {\em {Arch. Rational Mech. Anal.}}, 2015.

\bibitem{Fliss2021}
S.~Fliss, P.~Joly, and V.~Lescarret.
\newblock A {D}t{N} approach to the mathematical and numerical analysis in
  waveguides with periodic outlets at infinity.
\newblock {\em Pure Appl. Anal.}, 3(3):487--526, 2021.

\bibitem{Furuya2020}
T.~Furuya.
\newblock Scattering by the local perturbation of an open periodic waveguide in
  the half plane.
\newblock {\em J. Math. Anal. Appl.}, 489(1), 2020.

\bibitem{Hadda2016}
H.~Haddar and T.~P. Nguyen.
\newblock {A volume integral method for solving scattering problems from
  locally perturbed infinite periodic layers}.
\newblock {\em Appl. Anal.}, 96(1):130--158, 2016.

\bibitem{Hoang2011}
V.~Hoang.
\newblock The limiting absorption principle for a periodic semin-infinite
  waveguide.
\newblock {\em SIAM J. Appl. Math.}, 71(3):791--810, 2011.

\bibitem{Hohag2013}
T.~Hohage and S.~Soussi.
\newblock Riesz bases and {J}ordan form of the translation operator in
  semi-infinite periodic waveguides.
\newblock {\em J. Math. Pures Appl.}, 100(1):113--135, 2013.

\bibitem{Hu2021}
G.~Hu, W.~Lu, and A.~Rathsfeld.
\newblock Time-harmonic acoustic scattering from locally perturbed periodic
  curves.
\newblock {\em SIAM J. Appl. Math.}, 81(6):2569--2595, 2021.

\bibitem{Joly2006}
P.~Joly, J.-R. Li, and S.~Fliss.
\newblock Exact boundary conditions for periodic wave\-guides containing a
  local perturbation.
\newblock {\em Commun. Comput. Phys.}, 1:945--973, 2006.

\bibitem{Kirsc2019a}
A.~Kirsch.
\newblock Scattering by a periodic tube in $\mathbb{R}^3$: part i. the limiting
  absorption principle.
\newblock {\em Inverse Problems}, 35(10):104004, 2019.

\bibitem{Kirsc2019b}
A.~Kirsch.
\newblock Scattering by a periodic tube in $\mathbb{R}^3$: part i. the
  radiation condition.
\newblock {\em Inverse Problems}, 35(10):104005, 2019.

\bibitem{Kirsc2022}
A.~Kirsch.
\newblock A scattering problem for a local perturbation of an open periodic
  waveguide.
\newblock {\em Math. Methods Appl. Sci.}, 45(10):5737--5773, 2022.

\bibitem{Kirsc2023}
A.~Kirsch.
\newblock On the scattering of a plane wave by a perturbed open periodic
  waveguide.
\newblock {\em Math. Methods Appl. Sci.}, 46(9):10698--10718, 2023.

\bibitem{Kirsc2017}
A.~Kirsch and A.~Lechleiter.
\newblock The limiting absorption principle and a radiation condition for the
  scattering by a periodic layer.
\newblock {\em SIAM J. Math. Anal.}, 50(3):2536--2565, 2018.

\bibitem{Kirsc2017a}
A.~Kirsch and A.~Lechleiter.
\newblock A radiation condition arising from the limiting absorption principle
  for a closed full‐ or half‐waveguide problem.
\newblock {\em Math. Meth. Appl. Sci.}, 41(10):3955--3975, 2018.

\bibitem{Lechl2016}
A.~Lechleiter.
\newblock The {F}loquet-{B}loch transform and scattering from locally perturbed
  periodic surfaces.
\newblock {\em J. Math. Anal. Appl.}, 446(1):605--627, 2017.

\bibitem{Lechl2015e}
A.~Lechleiter and D.-L. Nguyen.
\newblock {Scattering of {H}erglotz waves from periodic structures and mapping
  properties of the {B}loch transform}.
\newblock {\em {Proc. Roy. Soc. Edinburgh Sect. A}}, 231:1283--1311, 2015.

\bibitem{Lechl2016a}
A.~Lechleiter and R.~Zhang.
\newblock A convergent numerical scheme for scattering of aperiodic waves from
  periodic surfaces based on the {F}loquet-{B}loch transform.
\newblock {\em SIAM J. Numer. Anal}, 55(2):713--736, 2017.

\bibitem{Lechl2017}
A.~Lechleiter and R.~Zhang.
\newblock A {F}loquet-{B}loch transform based numerical method for scattering
  from locally perturbed periodic surfaces.
\newblock {\em SIAM J. Sci. Comput.}, 39(5):B819--B839, 2017.

\bibitem{Yu2022}
X.~Yu, G.~Hu, W.~Lu, and A.~Rathsfeld.
\newblock {P}{M}{L} and high-accuracy boundary integral equation solver for
  wave scattering by a locally defected periodic surface.
\newblock {\em SIAM J. Numer. Anal.}, 60(5):2592--2625, 2022.

\bibitem{Yuan2007}
L.~Yuan and Y.~Y. Lu.
\newblock A recursive doubling {D}irichlet-to-{N}eumann map method for periodic
  waveguides.
\newblock {\em J. Lightwave technol.}, 25:3649--3656, 2007.

\bibitem{Zhang2017e}
R.~Zhang.
\newblock A high order numerical method for scattering from locally perturbed
  periodic surfaces.
\newblock {\em SIAM J. Sci. Comput.}, 40(4):A2286--A2314, 2018.

\bibitem{Zhang2019b}
R.~Zhang.
\newblock Numerical method for scattering problems in periodic waveguides.
\newblock {\em Numer. Math.}, 148:959--996, 2021.

\bibitem{Zhang2019a}
R.~Zhang.
\newblock Spectrum decomposition of translation operators in periodic
  waveguide.
\newblock {\em SIAM Journal on Applied Mathematics}, 81(1):233--257, 2021.

\bibitem{Zhang2021b}
R.~Zhang.
\newblock Exponential convergence of perfectly matched layers for scattering
  problems with periodic surfaces.
\newblock {\em SIAM J. Numer. Math.}, 60(2):804--823, 2022.

\bibitem{Zhang2021a}
R.~Zhang.
\newblock High order complex contour discretization methods to simulate
  scattering problems in locally perturbed periodic waveguides.
\newblock {\em SIAM J. Sci. Comput.}, 44(5):B1257--B1281, 2022.

\bibitem{Zhang2023}
R.~Zhang.
\newblock Higher-order convergence of perfectly matched layers in
  three-dimensional biperiodic surface scattering problems.
\newblock {\em SIAM J. Numer. Anal.}, 61(6):2917--2939, 2023.

\end{thebibliography}

\end{document}